\documentclass[11pt,a4paper]
{amsart}
\usepackage{enumerate, amssymb, paralist}
\usepackage[colorlinks]{hyperref}
\usepackage{hyperref}
\usepackage[dvipsnames]{xcolor}
\usepackage{graphicx}

\usepackage[dvipsnames]{xcolor}

\def\bdi{\begin{diagram}}
\def\edi{\end{diagram}}
\def\bsit{\begin{sit}}
\def\esit{\end{sit}}
\def\brem{\begin{rem}}
\def\brems{\begin{rems}}
\def\erem{\end{rem}}
\def\erems{\end{rems}}
\def\bprob{\begin{prob}}
\def\eprob{\end{prob}}
\def\bprobs{\begin{probs}}
\def\eprobs{\end{probs}}
\def\bques{\begin{ques}}
\def\eques{\end{ques}}
\def\bexa{\begin{exa}}
\def\bexas{\begin{exas}}
\def\eexa{\end{exa}}
\def\eexas{\end{exas}}
\def\bdefi{\begin{defi}}
\def\edefi{\end{defi}}
\def\bdefis{\begin{defis}}
\def\edefis{\end{defis}}
\def\bcor{\begin{cor}}
\def\ecor{\end{cor}}
\def\blem{\begin{lem}}
\def\elem{\end{lem}}
\def\bconv{\begin{conv}}
\def\econv{\end{conv}}
\def\bconj{\begin{conj}}
\def\econj{\end{conj}}
\def\bprop{\begin{prop}}
\def\eprop{\end{prop}}
\def\bthm{\begin{thm}}
\def\ethm{\end{thm}}
\def\bnota{\begin{nota}}
\def\enota{\end{nota}}
\def\bsit{\begin{sit}}
\def\esit{\end{sit}}
\def\be{\begin{equation}}
\def\ee{\end{equation}}
\def\bproof{\begin{proof}}
\def\eproof{\end{proof}}
\def\ba{\begin{array}}
\def\ea{\end{array}}
\def\la{\label}
\def\bcl{\begin{claim}}
\def\ecl{\end{claim}}

\newtheorem{thm}{Theorem}[section]
\newtheorem{cor}[thm]{Corollary}
\newtheorem{lem}[thm]{Lemma}
\newtheorem{prop}[thm]{Proposition}
\newtheorem{claim}[thm]{Claim}

\theoremstyle{definition}
\newtheorem{defi}[thm]{Definition}
\newtheorem{defis}[thm]{Definitions}
\newtheorem{conj}[thm]{Conjecture}

\newtheorem{conv}[thm]{Convention}
\newtheorem{nota}[thm]{Notation}
\newtheorem{rem}[thm]{Remark}
\newtheorem{rems}[thm]{Remarks}
\newtheorem{exa}[thm]{Example}
\newtheorem{exas}[thm]{Examples}

\newtheorem{sit}[thm]{}
\newcommand{\rien}[1]{}

\newcommand{\C}{\ensuremath{\mathbb{C}}}

\newcommand{\cO}{{\ensuremath{\mathcal{O}}}}

\newcommand{\Gr}{\rm Gr}

\def\NN{{\mathbb N}}
\def\ZZ{{\mathbb Z}}
\def\QQ{{\mathbb Q}}

\def\CC{{\mathbb C}}

\def\PP{{\mathbb P}}

\def\deg{\mathop{\rm deg}}

\def\Pic{\mathop{\rm Pic}}

\title[Lines, conics,  and all that]{
Lines, conics, and all that}

\author{C.~Ciliberto, 
M.~Zaidenberg}
\address{Dipartimento di Matematica, Universit\`a degli
Studi di Roma ``Tor Vergata'', Via della Ricerca Scientifica,
00133 Roma, Italy
}
\email{cilibert@mat.uniroma2.it}
\address{Univ.\ Grenoble Alpes, CNRS, Institut Fourier, F-38000 Grenoble, France}
\email{Mikhail.Zaidenberg@univ-grenoble-alpes.fr}

\dedicatory{To Bernard Shiffman on occasion of his seventy fifths birthday}

\thanks{
{\bf Acknowledgements:} This research started during a visit of the second author at the Department of Mathematics of the University Roma ``Tor Vergata'' in October, 2018. The authors acknowledge the MIUR Excellence Department Project awarded to the Department of Mathematics, University of Rome Tor Vergata, CUP E83C18000100006. 
Ciro Ciliberto is a member of GNSAGA of INdAM}

\thanks{\mbox{\hspace{11pt}}{\it 2010 Mathematics Subject Classification}:
14H45, 14N10, 14N25, 14J26, 14J28, 14Q05,  
14Q10.\\
\mbox{\hspace{11pt}}{\it Key words}: line, conic, rational curve, projective hypersurface, Fano threefold, Abel-Jacobi mapping}

\date{}
\begin{document}

\begin{abstract} This is a survey on the Fano schemes of linear spaces, conics, rational curves, and curves of higher genera in smooth projective hypersurfaces,   complete intersections, Fano threefolds, on the related Abel-Jacobi mappings, etc. 
\end{abstract}

\maketitle

\tableofcontents

\vfuzz=2pt
\thanks{}


\section*{Introduction} 
A classical problem of enumerative geometry asks to count curves with given numerical invariants in 
a smooth complete intersection variety $X$ in $\PP^n$. This includes the study of  various Fano schemes of $X$, that is, 
the components of the Hilbert schemes of curves with given numerical invariants, in particular, the Fano schemes of lines and conics. 
The present paper is a survey on this problem. We concentrate mainly on concrete numerical results. 
A special attention is paid
to the case of surfaces and threefolds. We discuss the lines and conics in Fano threefolds, 
which are not necessarily complete intersections. 
The latter  involves the Abel-Jacobi mapping, the related cylinder homomorphism, and  Torelli type theorems.  
To keep a reasonable volume, we do not outline proofs, and
restrict the exposition to smooth projective varieties. We are working over the complex number field unless otherwise is stated, 
while many results remain valid over any algebraically closed field of characteristic zero.
In the majority of cases, the names of the authors of cited results appear only in the list of references.
Together with the best known result, we  often mention its predecessors;
we apologize for possible confusions. 
We hope to provide the reader with a brief introduction to several chapters of this beautiful subject 
and a practical guide on the extended bibliography. 

Several important topics remain outside our survey,  such as, e.g.: rational curves in K\"{a}hler varieties \cite{CH}, in Moishezon varieties, see, e.g.,  \cite{Ogi}, in hyperk\"{a}hler varieties, see, e.g., \cite{AV}, in fibered varieties, see, e.g., \cite{An, An1}, pseudo-holomorphic curves in symplectic manifolds, see, e.g.,  \cite{AuLa, Gr, Hum, HT0, McDSa, McDSa1, Pansu}, minimal rational curves and the varieties of minimal rational tangents, see, e.g., \cite{Ar, ArCa, Hw, Hw1, HwKe, HwKi, HwMo1, Wa}, (minimal) rational and elliptic curves in flag varieties, Schubert varieties, Bott-Samuelson varieties, etc., see, e.g., \cite{BrKa, CoCo, HwMo2, LaMa, OSCW, PaPe, Pe0, Pe1, Wa1, Wa2}, varieties covered by lines and conic-connected varieties, see, e.g., \cite{IoRu1, IoRu2, IoRu3, MM, MMT}, Prym varieties, see, e.g., \cite{Bea0, Bea, De-3, DoSm, En, FS1, FS2, Iz, Mum, Rad1, Rama, To, We2, We3}, etc. We avoid the vast domain of Mirror Symmetry and the Gromov-Witten invariants, see, e.g., \cite{AsMo, Behr, CoKa, H1, H2, H3, H4T, Kon, KoMa, LiTe, MaPa, Prz}.  Among the topics missing in this survey, one should mention as well the results inspired by the Manin conjectures on the asymptotic of rational points and rational curves on Fano varieties \cite{Manin1, Manin2}, see, e.g., \cite{BaMa, BaTs, Bou0, Bou1, Bou2, KLO, LeTa, PaPe, Pe0, Pe1, Pey0, Pey1, Pey2, San, Saw, Tes05, Tes08, Tes09} and the literature therein. 

\section{Counting lines on surfaces} \la{ss:lines}  Cayley and Salmon \cite{Cal, Sa}, 
and also Clebsch \cite{Cl}, discovered that any smooth complex cubic surface in $\PP^3$ contains exactly 27 lines. See, e.g.,\cite{Hen} for a historical summary,  \cite{Dol, RaSt} for the modern treatment, \cite{BaSe, BLLP, Cal1, FiKh, KaWi, Manin0, OkTe, Se0, Sot}  for the count of lines on cubic surfaces over $\mathbb{R}$ and in characteristic $p>0$, \cite{HaKa} and \cite{McKMZ} for the monodromy of lines in families of cubic surfaces, and \cite{Cheng} for lines on singular cubic surfaces. 

In the case of quartic surfaces  in $\PP^3$ the following is known. 

\bthm\la{thm:quartics} \begin{itemize}\item[{\rm (a)}]  
The maximal number of lines on a smooth quartic in $\PP^3$ is $64$. This maximum is achieved by the F.~Schur quartic. 
Any smooth quartic with $64$ lines is  isomorphic to the F.~Schur quartic  {\rm \cite{DeItSe, RaSc3, Sc, Se1}}. 
\item[{\rm (b)}]
In the projective space $|\cO_{\PP^3} (4)|$ parameterizing the quartics in $\PP^3$, the subvariety of quartics containing a line has codimension one and degree $320$. The general point of this subvariety represents a quartic surface with a unique line  {\rm \cite{Man0, MaPa}}. 
\end{itemize}
\end{thm}

The first claim in (a) is due to 
B.~Segre \cite{Se1}, but his proof contains a mistake. 
In \cite{DeItSe} and \cite{RaSc3}  there are two different proofs based 
on the ideas of  B.~Segre. 

In the case of quintics we have the following weaker results. We let $\Sigma(d,n)$ stand for the projective space $|\cO_{\PP^n} (d)|=\PP^{{{n+d}\choose{n}}-1}$ parameterizing the degree $d$ 
effective divisors on $\PP^n$. The general point of $\Sigma(d,n)$ corresponds to a smooth hypersurface of degree $d$ in $\PP^n$.

\bthm\la{thm:quintics} 
 \begin{itemize}\item[{\rm (a)}] A smooth quintic surface in $\PP^3$ contains at most $127$ lines  {\rm \cite{RaSc5}}. 
 \item[{\rm (b)}]  The variety of quintic surfaces in $\PP^3$ 
 containing a line is irreducible of degree $1990$ and of codimension 
 $2$ in $\Sigma(5,3)$.  The general point of this variety corresponds to a quintic surface with a unique line {\rm \cite{MaRoXaVa, Man0}}. 
\end{itemize}
\ethm

The exact upper bound for the number of lines in a smooth quintic is unknown. 
The Fermat quintic and the Barth quintic contain exactly $75$ lines each.  

 For higher degree surfaces in $\PP^3$ the following is known (cf.\ Theorem \ref{thm:quartics}.b). 

\bthm\la{thm:surfaces-degree-d} 
\begin{itemize}\item[{\rm (a)}] 
A smooth surface of degree $d\ge 3$ in $\PP^3$ contains at most $11d^2 - 32d + 24$
 lines \cite{BauRa}.\footnote{This improves the former upper bound $11d^2-28d +12$ \cite{Se1}.}
\item[{\rm (b)}]
The surfaces of degree $d$ in $\PP^3$
containing a line are parameterized by an irreducible subvariety in $\Sigma(d,3)$ of codimension $d-3$ and degree
\be\la{eq:lines-on-surfaces} \frac{1}{24} {{d+1}\choose{4}} (3d^4 +6d^3 +17d^2 +22d+24)\,.\ee
The general such surface contains a unique line  {\rm \cite{MaRoXaVa, Man0}}. 
 \end{itemize}
 \ethm
 
 Once again, for $d\ge 5$ the exact upper bound in (a) is unknown. 
 The Fermat surface of degree $d$ contains exactly $3d^2$ lines. 
 As for (b), there are analogous formulas for the degrees of the loci of  degree $d$ surfaces in $\PP^3$ 
 passing through a given line and containing an extra line, and, respectively,  containing three lines in general position; 
 see \cite[Prop.\ 3.3 and 6.4.1]{MaRoXaVa}. The latter locus is  of codimension $3d+3-12=3d-9$ in $\Sigma(d,3)$. 
 For $d = 3$ its degree equals $720$, in agreement with the combinatorics of triplets of
skew lines contained in a  smooth cubic surface.
 
\section{The numerology of Fano schemes} \la{ss:lFano-scheme}   Given a variety $X\subset \PP^n$, the Fano scheme $F_h(X)$ 
is the scheme of $h$-planes contained in $X$.  
 Recall that a hypersurface is \emph{very general} if it belongs to the complement of countably many proper subvarieties in the space of hypersurfaces of a given degree. For the Fano schemes of lines on hypersurfaces one has the following results. 
 
 \bthm\la{thm:lines-in-hypersurfaces} 
 For a smooth hypersurface  $X$ of degree $d$ in $\PP^{n}$, where $n\ge 3$, the following holds. 
 \begin{itemize} 
\item[{\rm (a)}] {\rm  (\cite{BaVdV,  Beh1, Beh2, BehRi, LaRo, LaTo, Lew}, \cite[Ch.~V, 2.9, 4.3, 4.5]{Kol})}
\begin{itemize}
\item If $d\le 2n-3$ then the Fano scheme of lines $F_1(X) $ is nonempty;
\item if $d<n$ then $X$ is covered by lines. The latter holds also for non-smooth hypersurfaces;
\item $F_1(X)$ is smooth of the expected dimension 
$\delta=2n-3-d$ if either $X$ is general, or $d\le \min\{8,n\}$ and $\delta\ge 0$;
\item $F_1(X)$ is irreducible if $d\le (n+1)/2$  and $X$ is not a smooth quadric in $\PP^3$; it is connected if $d\le 2n-5$.
\end{itemize}
\item[{\rm (b)}] 
\begin{itemize}\item If  ${{d+1}\choose{2}}\le n$ then $F_1(X)$ is rationally connected and is a Fano variety for $X$ general \cite[Exerc.~V.4.7]{Kol}. 
\item If ${{d+2}\choose{2}}\ge 3n$ and $X$ is very general  then $F_1(X)$ contains no rational curve {\rm  \cite[Thm.~3.3]{RiYa2}}. 
\end{itemize}
\item[{\rm (c)}] Any hypersurface of degree $d=2n-3$ in $\PP^n$ contains a line. 
The number of lines in a general such hypersurface is \footnote{See also \cite{VdW} for another expression of this number.}
\[ d\cdot d!\sum_{k=0}^{n-3} \frac{(2k)!}{k!(k+1)!}\sum_{I\subset\{1,\ldots,n-2\},\,|I|=
n-2-k}\prod_{i\in I} \frac{(d-2i)^2}{i(d-i)}\,.\] 
For instance, a general quintic threefold in $\PP^4$ contains exactly $2875$ lines {\rm  \cite{Ha}}.
\end{itemize}
\ethm

The De Jong--Debarre Conjecture \cite{De3} states that $\dim(F_1(X))=\delta$ for any smooth hypersurface in $\PP^n$ of degree $d\le n$. 
The upper bound $n$ is optimal; indeed, for any $d > n$, there exists a smooth hypersurface $X$ of degree $d$ in 
$\PP^n$ with $\dim(F_1(X)) > \delta$ \cite[Cor.~3.2]{BehRi}. For instance, the Fano variety of lines on a Fermat hypersurface of degree $d\ge n$ in $\PP^n$ is of dimension $n-3$ \cite[Ch.~2, p.~50]{De0} (see also \cite{Ter}), which is larger than $\delta$ when $d > n$. On the Fermat quintic threefold $X_5$ in $\PP^4$ there are $50$ one-parameter
families of lines \cite{AlKa}, cf.\ (c). Whereas for $n=4$, $d=3$, the Fano variety of lines $F_1(X_3)$ on the cubic Fermat threefold $X_3=\{x^3_1 + x^3_2 + x^3_3 + x^3_4 + x^3_5 = 0\}$ in $\PP^4$ is a smooth surface (and so, $\dim(F_1(X_3))=\delta=2$) carrying 30 smooth elliptic Fermat cubic curves \cite{Ro0}.

For varieties with too many lines we have the following classification results, see  \cite[Thm.~4.2]{Cann}.

\bthm\label{thm:too-many-lines} {\rm (\cite{Rog, Se3})} 
Let $X\subset\PP^n$ be a $k$-dimensional subvariety of codimension at least $3$. Then $\dim(F_1(X))\le 2k -2$. Furthermore,
\begin{itemize}
\item if $\dim(F_1(X))=2k -2$ then $X = \PP^k$;
\item if $\dim(F_1(X))=2k -3$ then $X$ is either a quadric or a $1$-parameter family of $\PP^{k-1}$;
\item if $\dim(F_1(X))=2k -4$ then $X$ is either a $1$-parameter family of $(k - 1)$-dimensional quadrics, a $2$-parameter family of $\PP^{k-2}$, or the intersection of $6 - k$ hyperplanes with the Grassmannian ${\rm Gr}(2,5) \subset \PP^9$ in its Pl\"{u}cker embedding.
\end{itemize}
\end{thm}

For lines in an arbitrary subvariety in $\PP^n$ one has the following fact.

\bthm\label{thm:Landsberg}  {\rm (\cite[Thm.\ 1]{La}; see \cite{MePo} for $m=3)$} Let $X$ be a projective variety of dimension $m$ in $\PP^n$. Assume $X$ is covered by lines and through a general point of $X$ pass a finite number of lines. Then there are at most $m!$ lines passing through a general point of $X$.
\ethm

Notice that  a line osculating
to order $m+1$ at a general point of $X$ must be contained in $X$ \cite{La0}.  The bound $m!$ in the theorem is optimal; it is achieved for a smooth hypersurface of
degree $m$ in $\PP^{m+1}$. 
A similar uniform bound is known for higher degree curves. 

\bthm\label{thm:Hwang} {\rm (\cite[Thm.~2]{Hw0})}  Let $X$ be an irreducible projective variety of
dimension $m$ in $\PP^n$, and $x \in X$ be a general point. Let ${\rm Curves}_d(X, x)$ stand for the space of curves of degree $d$ lying on $X$ and passing through $x$. Then the number of components of ${\rm Curves}_d(X, x)$ is bounded above by
$${(2m+2)((2md)^m - 1)\choose {2m + 1}}^{(2m+2)(4d^2-4d+2)}\,\,.$$
\ethm

The exceptional locus in this theorem has codimension at least two in $X$ \cite[Thm.~3]{Hw0}.

For the Fano schemes of $h$-planes in hypersurfaces we have the following facts.

\bthm\label{thm:h-planes}  If $X$ is a smooth hypersurface of degree $d$ in $\PP^{n}$, where $n\ge 3$, then:
\begin{itemize}
\item $F_h(X)$ is irreducible of the expected dimension $\delta=(h+1)(n-h)-{{d+h}\choose{h}}$ provided \[2{{d+h-1}\choose{h}}\le n-h;\]
\item for $2h \ge \max\{(n-1)/2, n + 2 - d\}$ one has  {\rm \cite{Beh0}}
\[\dim (F_h(X)) \le \begin{cases} (m-h)(h+1) & \mbox{if}\,\,\, n=2m+1\\(m-h-1)(h+1) & \mbox{if}\,\,\,  n=2m
 \,;\end{cases}\]
\item
if $n$ is odd and $d>2$ then there are at most finitely many $(n-1)/2$-planes contained in $X$.  
If $n$ is even and $d>3$ then there is at most a one-parameter family of $(n/2-1)$-planes  in $X$ {\rm  \cite{Beh0}};
\item
assume $X$ is general, ${{d+h}\choose{h}}\le (n-h)(h+1)$,  and $d\neq 2$ or $n\ge 2h+1$. Then
there are explicit formulas for the degree of $F_h(X)$  {\rm \cite[Ex.~14.7.13]{DeMa},  \cite[Thm.~1.1]{H1}, \cite[Thms.~3.5.18, 4.3]{Man}}. \footnote{In the case of equality $F_h(X)$ is zero-dimensional and $\deg (F_h(X))$ is the number of $h$-planes in $X$.}
\end{itemize}
\ethm

If $X$ is a complete intersection  of type (that is, of multidegree) $\underline{d}=(d_1,\ldots,d_s)$ then the expected dimension of $F_h(X)$ is
\be\label{eq:ex-dim-ci}\delta=\delta(\underline{d},n,h)=(n-h)(h+1)-\sum_{i=1}^s {{d_i+h}\choose{h}}\,.\ee
 In particular, the expected dimension of $F_1(X)$ is
 \[\delta= 2(n-1)-(d+s),\quad{where}\quad d=\sum_{i=1}^s d_i\,.\]
 Let also
\[ \delta_-=\min\{\delta, n-2h-s\}\,.\]
For (a) and (b) of the following theorem see \cite{Bor, DeMa, La, Pr} and the references therein.

\bthm\la{thm:Fano-scheme} For a complete intersection $X\subset\PP^n$  of type $\underline{d}=(d_1,\ldots,d_s)$ the following holds. 
\begin{itemize}
\item[{\rm (a)}]  If $\delta_-<0$ then $F_h(X)=\emptyset$ for a general $X$.
\item[{\rm (b)}] If $\delta_-\ge 0$ then $F_h(X)$ has dimension $\delta$, 
is smooth for a general $X$, and is irreducible if $\delta_->0$. 
\item[{\rm (c)}]  If $\delta_-\ge 0$ and $\delta\ge n-h-s$ then  through any point $x \in X$ passes an $h$-plane contained  in $X$ {\rm \cite{Miz}}.
\item[{\rm (d)}] For any smooth Fano complete intersection $X\subset\PP^n$  of type $\underline{d}=(d_1,\ldots,d_s)$ with $s\ge 2$ and $d=\sum_{i=1}^s d_i\le s+5$, $F_1(X)$ has the expected dimension $\delta=2n-d-s-2$ {\rm \cite{Cann}}.
\end{itemize}
 \ethm
 
 See also \cite{ChIl, Ilten, IlKe, IlSu, IlZo} for the Fano schemes of toric varieties and of complete intersections in toric varieties, and \cite{Lar, LaKi} for applications to the machine learning.
 
\bthm\label{prop:ci}  {\rm (\cite{BCFS, Miz}; see \cite{Man} for $s=1$)}  Let $\Sigma(\underline{d},n)$ be the scheme which parameterizes the 
complete intersections of type $\underline{d}=(d_1,\ldots,d_s)$ in $\PP^n$, and let $\Sigma(\underline{d},n,h)$ be the subvariety of $\Sigma(\underline{d},n)$ of points which correspond to the complete intersections which carry $h$-planes. Then 
\[\dim \Sigma(\underline{d},n) = \sum_{j=1}^{s} {{d_j+n}\choose{d_j}}\,.\] If $\gamma(\underline{d},n,h):=-\delta(\underline{d},n,h)>0$ then 
 $\Sigma(\underline{d},n,h)$ is a nonempty, irreducible and rational subvariety of codimension  $\gamma(\underline{d},n,h)$ in $\Sigma(\underline{d},n)$. The general point of $\Sigma(\underline{d},n,h)$ corresponds to a complete intersection which contains a unique linear subspace of dimension $h$ and has singular locus of dimension $\max\{-1, 2h+s-n-1\}$ along its unique $h$-dimensional linear subspace (in particular, it is smooth provided $n\ge 2h+s$).
\ethm

To determine the degree of $\Sigma(\underline{d},n,h)$  we propose the following receipt. 

\bthm\label{prop:cig} {\rm (\cite[Cor.\ 2.2]{CZ})} Let $X$ be a general complete intersection of type $(d_1,\ldots, d_{s-1})$ in $\PP^n$ verifying
\[
\dim(F_h(X))>0\quad\text{and}\quad \dim(\Sigma(d_m,X))> \gamma(\underline{d},n,h)\geqslant 0\,,\]
where 
\[
\Sigma(d_m,X)=|\mathcal O_X(d_m)|\,.
\]
Let $\Sigma(d_s,X,h)$ be the set of points in $\Sigma(d_s,X)$ which correspond to complete intersections of type $\underline{d}=(d_1,\ldots,d_s)$ contained in $X$ and containing a linear subspace of dimension $h$.
Then the degree
$\deg(\Sigma(d_m,X,h))$ equals the coefficient of the monomial $x_0^ nx_1^ {n-1}\cdots x_h^ {n-h}$ in the product of the following polynomials in $x_0,\ldots, x_h$:
\begin{itemize}
\item the product $Q_{h,\underline{d}}=\prod_{i=1}^{s-1} Q_{h,d_i}$ of the polynomials
\[
Q_{h,d_i}=\prod_{v_0+\ldots+v_h=d_i} (v_0x_0+\cdots+v_hx_h)\,;
\]
\item the homogeneous component of degree 
\[
\rho:={{d_s+h}\choose h}-\gamma(\underline{d},n,h)=\dim(F_h(X))
\]
of the polynomial 
\[\prod_{v_0+\ldots+v_h=d_m}(1+v_0x_0+\ldots+v_hx_h)\,;\]
\item the Vandermonde polynomial $V(x_0,\ldots, x_h)$.
\end{itemize}
The general point of $\Sigma(d_m,X,h)$ corresponds to a complete intersection  of type $\underline{d}=(d_1,\ldots,d_s)$ which contains a unique subspace of dimension $h$.
\ethm

Notice that  our assumptions hold automatically if $\gamma(\underline{d},n,h)$ is sufficiently small, e.g., if $\gamma(\underline{d},n,h)=1$. 
 
\section{Geometry of the Fano scheme}

In this section we consider the complete intersections of type $\underline{d}=(d_1,\ldots,d_s)$ in $\PP^n$ whose Fano schemes have positive expected dimension $\delta=\delta(\underline{d}, n, h)>0,$ see \eqref{eq:ex-dim-ci}.
 We assume $d_i\ge 2$, $i=1,\ldots,s$. 
If also $n\ge 2h+s+1$ then by Theorem \ref{prop:ci}(b), for a general complete intersection $X$ of type $\underline{d}$ in $\mathbb P^ n$, the Fano variety $F_h(X)$ of linear subspaces of dimension $h$ contained in $X$ is a smooth, irreducible variety of dimension $\delta(\underline{d}, n, h)$. 

\bthm\label{prop:DM-formula} {\rm (\cite [Thm.\ 4.3]{DeMa})} In the notation and assumptions as before, 
 the degree of the Fano scheme $F_h(X)$ under the Pl\"ucker embedding equals the coefficient 
of the monomial $x_0^nx_1^{n-1}\cdots x_h^{n-h}$ of the product $Q_{h,\underline{d}}\cdot e^{\delta}\cdot V$ where 
\begin{itemize}
\item
$V$ stands for the Vandermonde polynomial $V({\bf x})=\prod_{0\le i< j \le h} (x_i-x_j)$;
\item $e({\bf x}):=x_0+\cdots+x_h$ and $\delta=\delta(\underline{d}, n, h)$;
\item 
 $Q_{h,\underline{d}}$ is  the product $\prod_{i=1}^{s} Q_{h,d_i}$ of the polynomials
\[
Q_{h,d_i}=\prod_{v_0+\ldots+v_h=d_i} (v_0x_0+\cdots+v_hx_h)\,.
\]
\end{itemize}
  \ethm
  
\brem\label{rem:Hiep} An alternative expression for $\deg(F_h(X))$  based on the Bott residue formula can be found in \cite[Formula (4)]{CZ}, \cite[Thm.\ 1.1]{H0}, and \cite [Thm.\ 2]{H2}; cf.\ also \cite[Ex.\ 14.7.13]{Ful}, \cite{Man0}, and \cite[Sect.\ 3.5]{Man}. 
\erem

There are formulas expressing certain numerical invariants of $F_h(X)$ other than the degree. If $\delta(\underline{d}, n, h)=1$ then $F_h(X)$ is a smooth curve; its genus was computed in \cite {H2}. In the case where $F_h(X)$ is a surface, that is, $\delta(\underline{d}, n, h)=2$, the Chern numbers of this surface and its holomorphic Euler characteristic $\chi(\mathcal{O}_{F_h(X)})$ were computed in \cite{CZ}. Actually, \cite{CZ} contains formulas for $c_1(F_h(X))$ and $c_2(F_h(X))$ in the general case provided $\delta(\underline{d}, n, h)>0$. Applying these formulas to the case where the Fano scheme is a surface, one can deduce  the classically known values and new ones, as in the following examples.

\bexas\label{ex:numerics}
 \begin{itemize}\item
In the case of  the Fano surface $F=F_1(X)$ of lines  on the general cubic threefold $X$ in $\PP^4$ one has \cite{AlKl, Li}, \cite[Ex.~4.1]{CZ}
\[\deg(F)=45,\quad e(F)=c_2(F)=27,\quad c_1(F)^2=K_F^2=45,\quad\text{and}\quad \chi(\mathcal{O}_F)=6\,.\] 
\item
For the Fano surface $F=F_1(X)$  of lines on a general quintic fourfold $X$ in $\PP^5$ one gets \cite[Ex.~4.2]{CZ}
\[\deg(F)=6125,\quad e(F)=c_2(F)=309375,\quad c_1^2(F)=K_F^2=496125,\quad\text{and}\quad 
\chi(\mathcal O_F)=67125\,.\]
\item
For the Fano surface $F=F_1(X)$ of lines  on the intersection $X$ of two general quadrics in $\PP^5$ one has \cite[Ex.~4.3]{CZ}
\[\deg(F) =32,\quad e(F)=c_2(F)=0,\quad\text{and}\quad K_F^2=c_1(F)^2=0\,.\] 
In fact, $F$ is an abelian surface \cite{Re}.  
\item
Finally, for  the Fano scheme $F=F_2(X)$ of planes on a general cubic fivefold $X$ in $\PP^6$ one gets \cite[Ex.~4.4]{CZ}
\[\deg(F)=2835,\quad e(F)=c_2(F)=1304,\quad  K_F^2=c_1(F)^2=25515,\quad\text{and}\quad \chi(\mathcal O_F)=3213\,.\]
\end{itemize}
\eexas

As for the Picard numbers of the Fano schemes of complete intersections, 
one has the following result. 

\bthm {\rm (\cite[Thm.\ 03]{ZJ}; cf.\ also \cite{DeMa})} Let $X$ be a very general complete intersection of type $\underline{d}=(d_1,\ldots,d_s)$ in $\PP^n$. Assume
$\delta(\underline{d}, n,  h) \ge 2$. 
Then $\rho(F_h(X)) = 1$ except in the following cases:
\begin{itemize}\item
$X$ is a quadric in $\PP^{2h+1}$, $h \ge 1$. Then $F_h(X)$ consists of two  isomorphic smooth disjoint components, and the Picard number of each component is $1$;
\item $X$ is a quadric in $\PP^{2h+3}$, $h \ge 1$. Then $\rho(F_h(X))=2$; 
\item $X$ is a complete intersection of two quadrics in $\PP^{2h+4}$, $h \ge 1$. Then $\rho(F_h(X))=2h + 6$.
\end{itemize}
\ethm

The assumption `very general' of this theorem cannot be replaced by `general'; one can find corresponding examples in \cite{ZJ}. See also \cite[Chapter~22]{Ha0} for the Fano varieties of smooth quadrics. It is well known, for instance, that the Fano variety of lines on a smooth three-dimensional quadric is isomorphic to $\PP^3$ \cite[Exercice~22.6]{Ha0}.\\

Next we turn to the irregular Fano schemes of general complete intersections. 

\begin{thm}\label{thm:irregular} {\rm (\cite[Thm.\ 5.1]{CZ})} Let $X$ be a general complete intersection of type $\underline{d}=(d_1,\ldots,d_s)$ in $\PP^n$. 
Suppose that the Fano scheme $F=F_h(X)$ of $h$-planes in $X$, $h\ge 1$, is irreducible of dimension $\delta\ge 2$. Then $F$ is irregular 
if and only if one of the following holds:
\begin{itemize}
\item[{\rm (i)}] $F=F_1(X)$ is the variety of lines on a general cubic threefold $X$ in $\PP^4$ {\rm ($\dim(F)=2$)};
\item[{\rm (ii)}] $F=F_2(X)$ is the variety of planes on a general cubic fivefold $X$ in $\PP^6$  {\rm ($\dim(F)=2$)};
\item[{\rm (iii)}] $F=F_h(X)$ is the variety of $h$-planes on the intersection  $X$  of two general quadrics in $\PP^{2h+3}$, $h\in\NN$  {\rm ($\dim(F)=h+1$)}. 
\end{itemize}
\end{thm}

\brem\label{rem:even-dim} 
1. 
The Fano surface of lines $F=F_1(X)$ on a smooth cubic threefold $X\subset \PP^4$ in (i) was studied by Fano \cite{Fa1} who found, in particular, that $q(F)=5$. Using Example \ref{ex:numerics} we deduce $p_g(F)=10$; cf.\ also \cite[Thm.\ 4]{Bea},  \cite{BSD, CG, Gh, 
Li, Ru, T1, T2}, \cite[Sect.~4.3]{Re}. There is an isomorphism ${\rm Alb}(F) \simeq J(X)$ where $J(X)$ is the intermediate Jacobian, see \cite{CG} and Section \ref{sec:q-lines}. The latter holds as well for $F=F_2(X)$ where $X\subset \PP^6$ is a smooth cubic fivefold as in (ii) \cite{Co}, cf. Example \ref{ex:numerics}. Thus, $q(F)>0$ in (i) and (ii). 

Consider further the Fano scheme $F=F_h(X)$ of $h$-planes on a smooth intersection $X$ of two quadrics in $\PP^{2h+3}$ as in (iii). By a  theorem of M.~Reid \cite[Thm.\ 4.8]{Re} (see also \cite[Thm.\ 2]{DR}, \cite{T4}), $F$ is isomorphic to the Jacobian $J(C)$ of a hyperelliptic curve $C$ of genus $g(C)=h+1$ (of an elliptic curve if $h=0$). Hence, one has $q(F)=\dim(F)=h+1>0$ for $h\ge 0$. There is an isomorphism $F\simeq  J(X)$ where $J(X)$ is the intermediate Jacobian of $X$ \cite{Do}.

2. The complete intersections in (i)-(iii) are Fano varieties. The ones in (i) are  Fano threefolds of index 
$2$ with a very ample generator of the Picard group. The Fano threefolds of index $1$ with 
a very ample anticanonical divisor which are complete intersections are the varieties $V^3_{2g-2}\subset \PP^{g+1}$ of genera  $g=3, 4, 5$, that is (see \cite[Ch. IV, Prop.\ 1.4]{Is}; cf.\  Section \ref{sec:q-lines}): 
\begin{itemize}
\item[$g=3$:] the smooth quartics $V^3_4$ in $\PP^4$;
\item [$g=4$:] 
the smooth intersections $V^3_6$ of a quadric and a cubic  in $\PP^5$;
\item[$g=5$:] the smooth  intersections $V^3_8$ of three quadrics in $\PP^6$. 
\end{itemize}
The Fano scheme of lines $F=F_1$ of a general Fano threefold 
$V^3_{2g-2}$ is a smooth curve of a positive genus $g(F)>0$. 
In fact, $g(F)=801$ for $g=3$, $g(F)=271$ for $g=4$, and $g(F)=129$ for $g=5$ \cite{Mar}, \cite[Ex.~1-3]{H2}, \cite[Thm.\ 4.2.7]{IP}. 
For $X=V^3_{2g-2}$ with $g\in\{3,4,5\}$, the Abel-Jacobi map 
$J(F) \to J(X)$ to the intermediate Jacobian is an epimorphism, 
see \cite{Is} and \cite[Lect.\ 4, Sect.\ 1, Ex.\ 1 and Sect.\ 3]{T3}.

3. Besides the Fano threefolds $V^3_{2g-2}$, there are other complete intersections whose Fano scheme of lines is a curve of positive genus. This holds, for instance, for the general hypersurface of degree $2r-4$ in $\PP^r$, $r\ge 4$, and for the general complete intersections of types $\underline{d}=(n-3,n-2)$ and $\underline{d}=(n-4,n-4)$  in $\PP^n$ where $n\ge 5$ and $n\ge 6$, respectively \cite[Ex.~1-3]{H2}. One can find in \cite{H2} a  formula for the genus of $F$ in these cases. 

4. Let $X$ be a smooth intersection of two quadrics in $\PP^{2k+2}$. Then the Fano scheme $F_k(X)$ is reduced and finite of cardinality $2^{2k+2}$ \cite[Ch.\ 2]{Re}, whereas $F_{k-1}(X)$ is a rational Fano variety of dimension $2k$ and index $1$, whose Picard number is $\rho = 2k + 4$, see \cite{AC, Ca}, and the references therein.  

5. It is shown in \cite{CCN, Ram} that the component of the Hilbert scheme whose general point is a pair of transversal  linear subspaces  in $\PP^n$ of given dimensions $a,b$, is a Mori dream space; its effective and nef cones are described,  and it is determined as to when such a component is Fano.
\erem

\section{Counting conics in complete intersection}\la{ss:conics} 

A \emph{conic} in $\PP^n$ is a curve whose Hilbert polynomial is $2t + 1$. Any conic $C$ is contained in a unique plane and is either a smooth (reduced) plane conic, or a pair of distinct lines, or a double line (see, e.g., \cite[Lem.~2.2.6]{DL}). Thus, a pair of skew lines does not fit in our terminology. Likewise, the general member of the Hilbert scheme ${\rm Hilb_3}(\PP^n)$ of curves  in $\PP^n$ whose Hilbert polynomial is $3t + 1$ is a twisted cubic.

Any smooth cubic surface $S$  in $\PP^3$ contains exactly 27 pencils of conics, and any smooth conic in $S$ belongs to a unique such pencil. Hence the variety of conics in $S$ is reducible and consists of $27$ $\PP^1$-components. 

By contrast, the number of conics on a general quartic surface in $\PP^3$ is finite. Furthermore, one has

\begin{thm}\label{thm:conics-in-quartics}  {\rm  (\cite{Bau})} There exist smooth quartic surfaces in $\PP^3$ which contain $432$ smooth conics; $16$ of these conics are mutually disjoint. 
\end{thm}

According  to \cite{Ni}, $16$ is the maximal number of disjoint rational curves on a quartic surface. 
In \cite{BaBa} and \cite{Ek} one can find constructions of two smooth quartic surfaces in $\PP^3$ carrying $352$ and $320$ smooth conics, respectively. 
The maximal number of conics  lying in a smooth quartic in $\PP^3$ is unknown. 
However, given any smooth quartic surface $S$, a general pencil of quartic surfaces through $S$ contains exactly $5016$ surfaces with a conic, counting things with multiplicities, see Theorem \ref{thm:conic-on-surfaces}. Let 
$\Sigma_c(d,3)$ be the variety of those degree $d$ surfaces in $\PP^3$ which contain conics. Then a general point of $\Sigma_c(4,3)$ corresponds to a smooth quartic surface carrying exactly two (coplanar) conics, and so, $\deg(\Sigma_c(4,3))=5016/2=2508$.

For higher degree surfaces in $\PP^3$ the following holds. 

\begin{thm}\label{thm:conic-on-surfaces} {\rm (\cite{Man0}; see also \cite[Prop.\ 7.1]{MaRoXaVa})} 
\begin{itemize}
\item
For $d\ge 4$, $\Sigma_c(d,3)$  is an irreducible subvariety of 
codimension $2d-7$ in $\Sigma(d,3)$.  In particular, 
$\Sigma_c(4,3)$ is a hypersurface of degree $2508$ in $\Sigma(4,3)$. \footnote{Cf.\ \cite{MaPa}.}
\item 
For $d \ge 5$, a general point of $\Sigma_c(d,3)$ corresponds to a surface which contains a unique (smooth) conic, and one has
 $$\begin{aligned} \deg(\Sigma_c(d,3))= {{d}\choose{4}} (d^2-d + 8)(d^2-d + 6) (207d^8-288d^7 + 498d^6 + 5068d^5 \\
-15693d^4 + 31732d^3-37332d^2 + 9280d-47040)/967680\,.\end{aligned} $$ 
\end{itemize}
\end{thm}

For higher dimensional hypersurfaces we have the following results.

\begin{thm}\label{thm:conic-on-hypersurfaces} \begin{itemize}\item[{\rm (a)}] Let $X$ be a hypersurface of degree
$n$ in $\PP^n$, $n\ge 2$. 
Then $X$ is covered by a family of conics. For the general such $X$,  the number of conics in $X$ 
passing through a general point equals 
\[\frac{(2n)!}{2^{n+1}} - \frac{(n!)^2}{2}\, \]  \cite[Prop.~3.2]{BoHo}, see also \cite[Thm.\ 0.1]{Lew},  
\cite{Te} $(n=4)$, \cite[Thm.\ 2]{CoMu} $(n=5)$. For $n\ge 4$ and a general $X$, the variety $R_2(X)$ of smooth conics in $X$ is smooth, irreducible, of dimension $n-2$. 
\item[{\rm (b)}] 
For the general hypersurface $X$ of degree $d$ in $\PP^n$, where $n\ge 3$, $R_2(X)$ 
is smooth of the expected dimension $\mu(d,n) := 3n-2d-2$ provided $\mu(d,n) \ge 0$, and is empty otherwise. 
It is irreducible provided $\mu(d,n)\ge 1$ and $X$ is not a smooth cubic surface in $\PP^3$ {\rm \cite[Thm.\ 1.1]{Fur1}}.
\item[{\rm (c)}]  If $X$  is general and $d \ge 3n - 1$, then $X$ contains no reducible conic \cite[Thm.~3.4]{RiYa2}.
\item[{\rm (d)}]  Let further $X$ be an arbitrary smooth hypersurface of degree $d$ in $\PP^n$.
\begin{itemize}\item 
Assume $n\ge 6$ and $d\le 6$ (so, $X$ is Fano). 
Let $R\neq\emptyset$ be an irreducible component of $R_2(X)$ such that the plane spanned by a general conic in $R$ is not contained in $X$. 
Then $\dim(R)=\mu(d,n)=3n-2d-2$  {\rm \cite{Fur2}}. 
\item Let $X$ be a smooth quartic threefold  in $\PP^4$ $(n=d=4)$. Then $\dim R\ge\mu(4,4)=2$
 for any irreducible component $R$ of $ R_2(X)$. Through a general point of $X$ passes $972$ conics  \cite{CoMuWe, Is-1, Te}.
\end{itemize}
\end{itemize}
\end{thm}

For instance, for a general sextic hypersurface $X$ in $\PP^5$ ($n=5$, $d=6$, and $X$ is Calabi-Yau) 
the Fano scheme $F_c(X)$ of conics in $X$ is a smooth projective curve, whose general point corresponds to a smooth conic 
\cite[Prop.~1.3-1.4]{Cao}. By (b),
for a general Fano hypersurface $X$ in $\PP^n$ (that is,  $d\le n$), $R_2(X)$ is smooth, irreducible, 
of dimension $\mu(d,n)\ge 1$ if $n\ge 4$. 

For conics in Fano complete intersections we have the following facts.

\begin{thm}\label{thm:var-conics-in-c.i.} {\rm (\cite[Lem.~1]{Shim1})} 
Let $X\subset\PP^{n}$ be the general complete intersection  of dimension $\geq 3$. 
Then the Fano scheme of conics in $X$ is a smooth projective variety. 
\ethm

\begin{thm}\label{thm:conics-in-c.i.} {\rm (\cite{Bea1, BoHo})}
Consider a smooth complete intersection $X\subset\PP^{n}$ 
of type $\underline{d}=(d_1,\ldots,d_s)$. 
Assume 
\begin{equation}\label{eq:relation} 2\sum_{i=1}^s d_i=n-s+1\,.
\end{equation}
Then the following holds.
 \begin{itemize}
\item
$X$ is a Fano variety with Picard number one of index $\iota(X)=\frac{n-s+1}{2}$. 
The anticanonical degree of a conic in $X$ equals $n-s+1$.  
\item
For a general $X$ the family of conics in $X$ is a nonempty, smooth and irreducible component 
of the Chow scheme of $X$ \footnote{The latter holds as well if the equality in \eqref{eq:relation}  
is replaced by the inequality ``$\le$'' \cite[Cor.~2.1]{BoHo}.}, of dimension $2(n-s-1)$. 
\item
Let $e(X)$
($e_0(X)$, respectively) be the number of conics passing through a general pair of points of $X$ 
(passing through a general point $x\in X$ and having a given general tangent direction at
$x$, respectively). Then for a general $X$ these conics are smooth, and one has
$$e_0(X) = e(X) = \prod_{i=1}^s (d_i-1)!d_i!\,.$$ 
\end{itemize} 
\end{thm}

See \cite{Bea1} for formulas for the numbers of lines and conics in $X$ 
meeting three general linear subspaces in $\PP^{r}$ of suitable dimensions; cf.\ also \cite{Li} 
in the case of lines,  \cite[Cor.\ 1.5]{BehKu}, and \cite{Va}.  See \cite{BGLP} for the variety of lines in $\PP^3$ tangent to four given quadric surfaces. See also \cite{Pan1} for a study on the Kontsevich moduli spaces $\overline{\mathcal{M}}_{0,2}(X, 2)$ of conics through a pair of points in a smooth complete intersection $X$ (see the definition below).

Set \[\epsilon(d,n)=2d+2-3n=-\mu(d,n)\,.\] Consider the subvariety $\Sigma_c(d,r)$ of $\Sigma(d,r)$ 
whose points correspond to hypersurfaces containing plane conics.

\begin{thm}\label{lem:tbpp}  {\rm (\cite[Thms.\ 6.1, 6.6]{CZ})} Assume $d\ge 2$, $n\ge 3$, and $\epsilon(d,n)\ge 0$. Then the following holds.
\begin{itemize}
\item[{\rm (a)}]  $\Sigma_c(d,n)$ is irreducible of codimension  $\epsilon(d,n)$ in $\Sigma(d,n)$. 
\item[{\rm (b)}]  If $\epsilon(d,n)> 0$ and $(d,n)\neq (4, 3)$ then the hypersurface corresponding to the general point of $\Sigma_c(d,n)$ contains a unique conic, and this conic is smooth.
In the case $(d,n)=(4, 3)$ it contains exactly two distinct conics, and these conics  are smooth and coplanar.
\item[{\rm (c)}] If $\epsilon(d,n)>0$ and $(d,n)\neq (4,3)$ then one has
\[\deg (\Sigma_c(d,n))= -\frac{5}{32}{{n+1}\choose{3}}\eta(1,1,1)\,,\]
where $\eta$ is the homogeneous form of degree $3n-1$ in the formal power series 
decomposition of 
\[
\left(\prod_{|{\bf v}|=d}(1+\langle {\bf v}, {\bf x}\rangle)\right)\\ \cdot  \left(\prod_{|{\bf v}|=d-2}(1+\langle {\bf v}, {\bf x}\rangle)\right)^ {-1}\,
\]
with ${\bf x}=(x_1,x_2,x_3)$, ${\bf v}=(v_1,v_2,v_3)\in (\ZZ_{\ge 0})^3$, and $|{\bf v}|=v_1+v_2+v_3$.
\end{itemize}
\end{thm}

The latter formulas are obtained by applying Bott's residue formula; see, e.g., \cite{Bo, Br, EG, H1, MV} for generalities. 

Recall the definition of the Kontsevich moduli spaces of stable maps \cite{Kon}. 
Let $X$ be a smooth projective variety in $\PP^n$. The \emph{Kontsevich moduli space} $\overline{\mathcal{M}}_{g,r}(X, e)$  parameterizes the isomorphism classes of corteges $(C, f, x_1,\cdots, x_r)$ where
\begin{itemize} \item
$C$ is a proper, connected, nodal curve of arithmetic genus $g$; 
\item $f\colon C\to X$ is a morphism whose image is a curve of degree $e$ in $\PP^n$;
\item $(x_1,\cdots, x_r)$ is an ordered collection of distinct smooth points of $C$; 
\item the cortege $(C, f, x_1,\cdots, x_r)$ admits only finitely many automorphisms.
\end{itemize}
In general, $\overline{\mathcal{M}}_{g,r}(X, e)$ is a proper Deligne-Mumford stack. The underlying variety of $\overline{\mathcal{M}}_{g,r}(X, e)$ is
projective, but does not  need to be smooth or irreducible. However, $\overline{\mathcal{M}}_{0,0}(X, e)$ is a compactification of the variety $R_e(X)$ of smooth rational curves in $X$ of degree $e$. 

Recall that a curve $C \simeq\PP^1$ in a smooth projective variety $X$ of dimension $m$
is called \emph{free} if $\mathcal{N}_{C/X} \simeq \mathcal{O}(a_1) \oplus \ldots \oplus \mathcal{O}(a_{m-1})$ where $a_i\ge 0$ $\forall i$. Statement (a) of the next theorem follows from Theorem \ref{thm:conic-on-hypersurfaces}(b).

\bthm\la{thm-DeLand} 
\begin{itemize}
\item[{\rm (a)}]
Let $X$ be  a smooth hypersurface  in $\PP^n$ of degree $d < n$. 
\begin{itemize}\item
 If $X$ is general then $\overline{\mathcal{M}}_{0,0}(X, 2)$ is irreducible and of the expected dimension $\mu(d,n)=3n-2d-2$  {\rm \cite[2.3.4]{DL}}. 
\item If $X$ is arbitrary (but smooth) then there is a unique component of $\overline{\mathcal{M}}_{0,0}(X, 2)$ which contains a conic passing through the general point of $X$. Moreover, if the dimension of the variety of non-free lines on $X$ is at most $n - 3$ then
there is a unique component of $R_2(X)$ whose general point corresponds to a smooth conic through the general point of $X$  {\rm \cite[2.3.6]{DL}}.
\end{itemize}
\item[{\rm (b)}]
Let $X$ be the Fermat hypersurface in $\PP^n$ of degree $d$. Then 
the moduli space $\overline{\mathcal{M}}_{0,0}(X, e)$ is irreducible of the expected dimension $e(n + 1 - d) + (n - 4)$ in the following cases:
\begin{itemize} 
\item
$e\ge 2$ and $ed \le n + 1$;
\item $e \ge 3$ and $ed \le n-1$.
\end{itemize}
It is irreducible as well if $e=3$ and $3d \le n+ 5$ {\rm \cite[2.4.4--2.4.6, 2.4.23, 2.4.30, 2.4.36]{DL}}. 
\end{itemize}
\ethm

\section{Lines and conics on Fano threefolds and the Abel-Jacobi mapping}\label{sec:q-lines}
\subsection{The Fano-Iskovskikh classification} The content of this section is partially borrowed from \cite[Sect.\ 4.1]{IP} and \cite[Sect.\ 2]{KuPrSh}. Let $X$ be a smooth Fano threefold, that is, a smooth threefold with an ample anticanonical divisor $-K_X$. One attributes to $X$ 
the following integers:
\begin{itemize}\item  the \emph{genus}\footnote{The genus $g(X)$
 is equal to the genus $g$ of a general curve section of the anticanonical model of $X$.}  \[g(X)= ( -K_X)^3/2 +1=\dim|-K_X|-1\ge 2\,;\] 
\item the \emph{index} $\iota(X)$, that is, the maximal natural number in $\{1,\ldots,4\}$ such that $-K_X=\iota(X)H$ for an ample divisor $H$ on $X$;
\item the \emph{Picard rank} $\rho(X)$ such that ${\rm Pic}(X)\simeq\ZZ^{\rho(X)}$;
\item the \emph{degree} $d(X)=H^3$;
\item the \emph{Matsusaka constant} $m_0=m_0(X)$, that is, the minimal integer such that $m_0H$ is very ample. 
\end{itemize}
One has $g(X)=\iota(X)^3d(X)/2+1$. Hence $d(X)=2g(X)-2$ if and only if $\iota(X)=1$.

The Fano-Iskovskikh classification of the Fano threefolds with Picard rank $\rho(X)=1$ yields 

\bthm\la{thm: classif-3folds} {\rm (\cite{Is}, \cite[Sect.\ 4.1]{IP})} 
Let $X$ be a smooth Fano threefold of genus $g$ with $\rho(X)=1$. Then one of the following {\rm (i)}--{\rm (iii)} holds.
\begin{itemize} \item[{\rm (i)}] The anticanonical divisor is very ample, and the linear system
$|-K_X|$ defines an embedding $\varphi$ of $X$ onto a projectively normal threefold $\varphi(X)$ of degree $2g - 2$ in $\PP^{g+1}$ with one of the following:
\begin{itemize}\item[$\rm (i3)$] $g=3$ and $\varphi(X)\subset\PP^4$ is a smooth quartic threefold; 
\item[$\rm (i4)$] $g=4$ and $\varphi(X)\subset\PP^5$ is complete intersection of a quadric and a cubic hypersurfaces; 
\item[$\rm (i5)$] $g=5$ and $\varphi(X)\subset\PP^6$ is complete intersection of three quadric  hypersurfaces;
 \item[$\rm (ig)$] $g\ge 6$ and $\varphi(X)\subset\PP^{g+1}$ is an intersection of quadric  hypersurfaces.
\end{itemize}
\item[{\rm (ii)}] $g=2$ and $X$ is a sextic double solid, that is, $|-K_X|$ defines a double cover $X \to \PP^3$ ramified along a smooth surface $S\subset\PP^3$ of degree $6$;
\item[{\rm (iii)}] $g=3$ and $|-K_X|$ defines a double cover $X \to Q$ over a smooth quadric threefold $Q\subset\PP^4$ ramified along a smooth surface $S\subset Q$ of degree $8$.
\end{itemize}
\ethm

The table of numerical data of  the Fano threefolds with $\rho(X)=1$ can be found in \cite[Sect.\ 12.2]{IP}. 
These threefolds form $17$ deformation families. According to the index, these are:
\begin{itemize}\item[$\iota=1$]: $10$ families with genera varying from $2$ to $12$ excluding $11$;
\item[$\iota=2$]: $5$ families of \emph{del Pezzo threefolds} with anticanonical degree $-K_X^3=8d$, $d=1,2,3,4,5$; 
\item[$\iota=3$]: the smooth quadric $Q$ in $\PP^4$ with anticanonical degree $54$; 
\item[$\iota=4$]: $\PP^3$  with anticanonical degree $64$. 
\end{itemize}

The families  of Fano threefolds $V^3_{2g-2}\subset\PP^{g+1}$ with $\rho=1$ and $\iota=1$  are classified according to the genus $g$ as follows. 
\begin{itemize}\item[$g=2,3,4,5$:] the Fano threefolds listed in (i3)-(i5), (ii), and (iii);
\item[$g=6$:] the smooth intersections ${\rm Gr}(2;5)\cap \PP^7\cap Q$ of the Grassmannian ${\rm Gr}(2;5)$
with a linear subspace $\PP^7$ and a quadric $Q$ in $\PP^9$, and 
\item[$g=6$:] the \emph{Gushel threefolds} of genus $6$, that is, the double covers of the del Pezzo threefold $Y={\rm Gr}(2;5)\cap \PP^6\subset\PP^9$ of degree $5$ branched along a smooth quadric section $Q\cap Y$;
\item[$g=7$:] the smooth linear sections ${\rm OGr}_+(5;10)\cap \PP^8$ of a connected component ${\rm OGr}_+(5;10)$ of the orthogonal Lagrangian Grassmannian in $\PP^{15}$;
\item[$g=8$:] the smooth linear sections ${\rm Gr}(2;6)\cap \PP^9$ of the Grassmannian ${\rm Gr}(2;6)\subset\PP^{14}$;   
\item[$g=9$:] the smooth linear sections ${\rm LGr}(3; 6)\cap \PP^{10}$  of the symplectic Lagrangian Grassmannian ${\rm LGr}(3; 6) \subset\PP^{13}$;
\item[$g=10$:] the smooth linear sections $\Omega^5\cap\PP^{11}$ of the homogeneous $G_2$-fivefold $\Omega^5\subset\PP^{13}$ 
(an adjoint orbit of the group $G_2$);
\item[$g=12$:] the smooth zero loci of triplets of sections of the rank $3$ vector bundle $\Lambda^2\mathcal{E}^{\vee}$, where $\mathcal{E}$ is the universal bundle over the Grassmannian ${\rm Gr}(3; 7)$.
\end{itemize}

Initially, the Fano threefolds $V^3_{2g-2}$ with $g=7,9,12$ were obtained from simpler ones via certain birational transformations (\emph{elementary Sarkisov links}); see, e.g., \cite[Prop.~3.4.1, 4.4.1, Thm.~4.3.3, 4.4.11]{IP}, \cite{Ta}.  We use above the Mukai description \cite{Mu0, Mu1, Mu2, Mu3} of the $V^3_{2g-2}$ with $g=7,\ldots,10$ as linear sections $X^n_{2g-2}\cap \PP^{g+1}$ of certain special Grassmannians $X^n_{2g-2}=G/P$, which are flag varieties  embedded in $\PP^{g+n-2}$.
 
Notice that the family of Gushel threefolds of genus $6$ (called also \emph{special Gushel-Mukai threefolds}) is a flat specialization of the family ${\rm Gr}(2;5)\cap \PP^7\cap Q$ of \emph{general Gushel-Mukai threefolds}. Thus, the Fano threefolds of genus $6$ form one deformation family of Gushel-Mukai threefolds. The same holds for the Fano threefolds of genus $3$; indeed, the family of smooth quartic threefolds in $\PP^4$ specializes to the double covers $X \to Q$ ramified along $Q\cap Y$, where $Y$ is a quartic  in $\PP^4$.  

 The families of Fano threefolds with $\rho=1$ and $\iota=2$ are classified according to the anticanonical degree $-K_X^3\in\{8d,\, d=1,2,3,4,5\}$ as follows (Fujita \cite{Fuj}; see \cite[Thm.~3.3.1]{IP}). 
\begin{itemize}
\item[$d=1$:] the smooth hypersurfaces of degree 6 in the weighted projective space $\PP(3,2,1,1,1)$. Another realization: the Veronese double cones, that is, the double covers 
$X\to V$, where $V\subset\PP^6$ is the cone over the second Veronese surface in $\PP^5$, branched at  the vertex $v$ of $V$ and along a smooth intersection of $V$  with a cubic hypersurface which does not pass through $v$ \cite{Is-1, Ti3}, see also \cite{Grin1, Grin2, HwKi};
\item[$d=2$:] the quartic double solids, that is, the double covers $X \to \PP^3$ branched along a smooth quartic surface $S\subset\PP^3$;
\item[$d=3$:] the smooth cubic threefolds $X \subset \PP^4$;
\item[$d=4$:] the smooth complete intersections of two quadrics in $\PP^5$;
\item[$d=5$:] the smooth linear sections ${\rm Gr}(2;5)\cap\PP^6$ of the Grassmannian ${\rm Gr}(2;5)\subset\PP^{9}$.
\end{itemize}

\subsection{Lines and conics on Fano threefolds} Let $X$ be a Fano variety of index $\iota(X)$, and let $H=K_X/\iota(X)\in\Pic(X)$. The \emph{lines} and \emph{conics} on $X$ are the irreducible curves $C$ in $X$ 
satisfying $C\cdot H=1$ and $C\cdot H=2$, respectively. One considers the Fano schemes $F_1(X)$  of lines and $F_c(X)$ of conics in $X$ 
meaning actually the unions of the components of the Hilbert schemes of $X$ whose general points correspond to lines and conics on $X$, respectively.

In the case where $-K_X$ is very ample, e.g., if $\rho(X)=\iota(X)=1$, the lines and conics on $X$ are sent to the usual lines and conics  under the anticanonical embedding $X\hookrightarrow\PP^{2g(X)-2}$. Otherwise, consider, for instance, a \emph{double solid}, that is, a double cover $\pi\colon X\to\PP^3$ branched along a smooth surface $S\subset\PP^3$ of degree $4$ 
 (degree $6$, respectively) which does not contain any line (any conic, respectively). An irreducible curve $C$  on $X$ is a line  (a conic, respectively) 
 if and only if $C'=\pi(C)$ is a bitangent line 
 of $S$ in $\PP^3$ (a conic in $\PP^3$ with only even local intersection indices with $S$, respectively). In these cases $\pi^*(C') = C + i(C)$ has two 
 irreducible components, where $i$ is the
involution associated to $\pi$. 

The following theorem summarizes results from \cite{CeVe, CoMuWe, Co, Ili94, Is, Let, Lo2, Sh, Ta, Te}.

\bthm\label{thm:lines-in-3folds}  
Let $X$ be a Fano threefold with $\rho(X)=1$ of index $\iota(X)=1$ and genus $g=g(X)$. Then the following holds.
\begin{itemize}\item Every line on $X$ meets $l(g)$ lines counting things with multiplicities, where 
\begin{itemize} 
\item $l(2)=625$;
\item $l(3)=81$ if $X$ of genus $g=3$ is a double cover of a quadric $Q\subset\PP^4$; 
\item $l(4)=31$;
\item $l(5)=17$; 
\item $l(6)=11$; 
\item $l(8)=6$  provided $X$ of genus $g=8$ is general  {\rm \cite{BlMu, Mar}}.
\end{itemize}
\item If $-K_X$ is very ample 
and $g\ge 3$ then $F_1(X)$ is of pure dimension $1$; it is smooth, reduced and irreducible for a  general $X$. 
Through a point of $X$ passes at most a finite number of lines if $g=3$ and $X$ is a general quartic surface in $\PP^4$, 
at most $6$ lines if $g=4$, and at most $4$ lines if $g\ge 5$ \cite[4.2.2, 4.2.7]{IP}.
\item Assume $-K_X$ is very ample and $g\ge 5$. Then $F_c(X)$ is two-dimensional. 
Furthermore, through almost any point of $X$ (any point if $g\ge 10$)  
passes at most a finite number of conics. 
A general conic in $X$ meets at most a finite number of lines if $g\ge 5$; the latter is true for any conic if $g\ge 9$ \cite[4.2.5--4.2.6]{IP}.
\end{itemize}
\end{thm}

Notice that there are smooth quartic threefolds in $\PP^4$ (for instance, the Fermat quartic) which contain cones over curves. However, through a general point of any smooth quartic in $\PP^4$ passes exactly $972$ conics \cite{CoMuWe, Is-1}. See also \cite[4.2.7]{IP} for the genera of the curve $F_1(X)$ for a general $X$ as in the theorem.  

The following two theorems summarize the results of  \cite{AlKl,  BF, DR, FN,  Ili03, IM1, Is, KuPrSh, Mar, Put, Te}.

\bthm\label{thm:q-lines} {\rm  (\cite[Thm.~1.1.1]{KuPrSh})} Let $X$ be a smooth Fano threefold with $\rho(X)=1$ of index $\iota=2$
and degree $d=d(X)\ge 3$. Then the Fano scheme of lines $F_1(X)$  is a smooth irreducible
surface. In particular, 
\begin{itemize}\item[$d= 3$:] $F_1(X)$  is a minimal surface of general type with irregularity $5$, geometric genus $10$, and canonical degree $K_{F_1(X)}^2=45$;
\item[$d=4$:] $F_1(X)$ is an abelian surface;
\item[$d=5$:] $F_1(X)\simeq\PP^2$. 
\end{itemize}
\ethm

\bthm\label{thm:q-conics}  {\rm  (\cite[Thm.\ 1.1.1]{KuPrSh})} 
Let $X$ be a smooth Fano threefold with $\rho(X)=1$  of index $1$ and genus $g=g(X) \ge 7$.
Then the Fano scheme of conics $F_c(X)$ is a smooth irreducible surface. More precisely,
\begin{itemize}\item[$g=7$:] $F_c(X)$ is  symmetric square of a smooth curve of genus $7$;
\item[$g=8$:] $F_c(X)$ is a minimal surface of general type with irregularity $5$, geometric genus $10$, and canonical degree $K_{F_c(X)}^2=45$;
\item[$g=9$:] $F_c(X)$ is a ruled surface isomorphic to the projectivization of a simple rank $2$ vector bundle on a smooth curve of genus $3$;
\item[$g=10$]: $F_c(X)$ is an abelian surface;
\item[$g=12$:] $F_c(X)\simeq \PP^2$.
\end{itemize}
\ethm

There exists the following duality between the Fano threefolds of indices 1 and 2 based on the Mukai construction \cite{Mu1, Mu2}. 

\bthm\label{thm:duality} {\rm (\cite{Kuz}; see also \cite[Appendix B]{KuPrSh})} For any smooth Fano threefold $X$ with $\rho(X)=1$, 
$\iota(X)=1$, and $g(X)\in\{8, 10, 12\}$ there is a smooth Fano threefold $Y$ with
$\rho(Y)=1$, $\iota(Y)=2$, and $d(Y)=g(X)/2-1$ such that $F_1(X)\simeq F_c(Y)$ (and the derived categories of $X$ and $Y$ are equivalent).
\ethm

Notice \cite{IlSc} that there is no nonconstant morphism from a Fano threefold with $\rho=1$ of index $1$ to a Fano threefold with $\rho=1$ of index $2$. This was conjectured by Th.~Peternell \cite{Pet}; the proof exploits the Fano schemes of lines and conics on these threefolds, in particular, the subfamilies of reducible conics, and the families of lines on a tangent scroll to a curve contained in these threefolds. \\

The following theorem covers some results of \cite{ChSh, Fae, FN, Ili94, San, TaZu}.

\bthm\label{thm:sections-of-Gr-2-5} {\rm (\cite{ChHoLe})} Consider the quintic Fano variety  $Y^m_5$ of dimension $m$, where $2\le m\le 6$, and of index $m-1$, which is the general linear section of the Grassmannian $\Gr(2,5)$ under its Pl\"{u}cker embedding in  $\PP^9$. Then  for any $d\in\{1,2,3\}$ the moduli space $R_d(Y_5^m)$ of smooth rational curves of degree $d$ on $Y_5^m$ is  rational. For $m = 4, 5$ the Hilbert scheme  ${\rm Hilb}_2 (Y_5^m) $ of conics on $Y_5^m$ is smooth and irreducible.
In particular,
 \[{\rm Hilb}_1(Y_5^3)= \PP^2, \,\,\, {\rm Hilb}_2(Y_5^3) = \PP^4,\,\,\,  {\rm Hilb}_3(Y_5^3)= \Gr(2, 5),\,\,\,\text{and}\,\,\, {\rm Hilb}_1(Y_5^6)=R_1(Y_5^6) =  \Gr(1, 3, 5)\]  is a flag variety.
 \ethm

For the Fano schemes of lines and conics on del Pezzo threefolds, Fano fourfolds, and some other higher dimensional Fano varieties see, e.g., \cite{BD, DIM2, DoMa, Do, HT1, KaRa, MaSt, MT3, PrZa}; see also \cite{Suz} for polarized Fano varieties covered by linear spaces.\\

Let $G$ be simple linear algebraic group, $P\subset G$ be a parabolic subgroup. Consider the flag variety $X=G/P$ along with its minimal homogeneous embedding in a projective space. Then the Fano schemes $F_k(X)$ of linear subspaces in $X$ are disjoint unions of flag varieties \cite[Sec.~4.2--4.3]{CoCo, LaMa}; see also \cite{Man1}.

\subsection{The Abel-Jacobi mapping}\la{ss: interm-Jac}  Recall \cite{CG} that the \emph{intermediate Jacobian} \[J(X) = H^{2,1} (X)^*/(H_3(X, \ZZ)\,\text{modulo torsion})\] of a smooth Fano threefold $X$ is a principally polarized abelian variety. Using a fine structure of the intermediate Jacobian one can detect the non-rationality of $X$; see, e.g., \cite{Bea4, CG, GaSh, Kol2, Pro2, Voi4, Wi}, cf. also \cite{AM, ChP, Pro1, Ro}. Given a variety $F$ we let $A(F)$ be the Albanese variety of $F$. For a Fano scheme $F(X)$ of a Fano threefold $X$ the \emph{Abel-Jacobi mapping} $A(F(X))\to J(X)$ is defined via the \emph{cylinder homomorphism} $H_1(F(X),\ZZ)\to H_3(X, \ZZ)$, see the next section. For the Fano threefolds $X$  one considers the Fano surface of lines $F(X)=F_1(X)$ if $\iota(X)=2$ and of conics $F(X)=F_c(X)$ if $\iota(X)=1$. 
For certain Fano threefolds $X$, 
the Abel-Jacobi mapping is known to be either an isomorphism, or an isogeny, hence $q(F)=h^{2,1}(X)$. We summarize these results in the following theorem. See Table 12.2 in \cite{IP} for the values of $h^{2,1}(X)$.

\bthm\label{thm-Fc} {\rm (\cite[\S 8.2]{IP})} Consider a Fano threefold  $X$ with $\rho(X)=1$ of genus $g=g(X)$, degree $d=d(X)$, and  index $\iota=\iota(X)\in\{1,2\}$. 
\begin{itemize}\item[$\iota=2$:] Assume $\iota=2$, and  let $F=F_1(X)$ be the Fano scheme of lines on $X$. Then the following holds.
\begin{itemize} \item[$d=5$:] $A(F)$ and $J(X)$ are both trivial \cite{Is}.
\item[$d=2,3,4$:] $F$ is a smooth irreducible surface, the Abel-Jacobi mapping $A(F)\to J(X)$ is an isomorphism, and $X$ is uniquely determined by $F$ in the following cases: 
\begin{itemize} 
\item[$d=2$:] 
$X\to\PP^3$  is a quartic double solid whose branching surface $S\subset\PP^3$ has no line. One has $q(F)=10$ 
{\rm \cite[Ch. III, Sect.\ 1]{Is}, \cite{Clm0, Clm5, Ti2, Ti1,  We, We1}}; 
\item[$d=3$:]   $X$ is a smooth cubic threefold in $\PP^4$. One has $q(F)=5$ {\rm \cite{CG}};
\item[$d=4$:] $X$ is a smooth complete intersection of two quadrics in $\PP^5$. One has $q(F)=2$, $F\cong J(X)$ is an abelian surface \cite{DR, Do, Re, T4}. 
\end{itemize}
\item[$d=1$:] Let  $X\to V$ be a double Veronese cone branched along a smooth surface $W\subset V$ 
cut out by a cubic hypersurface in $\PP^6$, and let $F_0$ be the Hilbert scheme of conics in $V$ $3$-tangent to $W$. Then $F$ and $F_0$ are  smooth irreducible surfaces, there is a branched double covering $\pi\colon F\to F_0$,  the Fano scheme $\mathcal{F}(X)$ is not reduced and consists of $F$ and the embedded ramification curve of $\pi$, and the Abel-Jacobi mapping yields an isogeny $A(F)/\pi^*A(F_0)\to J(X)$ where $\dim J(X)=h^{2,1}(X)=21$, and $X\to V$  is uniquely determined by the pair $(F,\pi)$ \cite{Ti3}.
\end{itemize}
\item[$\iota=1$:] Assume $\iota=1$, and let $F=F_c(X)$  be the Fano scheme of conics on $X$. Then the following holds.
\begin{itemize} 
\item[$g=2,3,4,5$:] $F$ is a smooth irreducible surface and the Abel-Jacobi mapping $A(F)\to J(X)$ is an isomorphism in the following cases: 
\begin{itemize} 
\item[$g=2$:]  $X\to\PP^3$
is the general sextic double solid. One has $q(F)=52$ {\rm \cite[Thm.\ 3.3]{CeVe}};
\item[$g=3$:]  $X$ is the general quartic threefold  in $\PP^4$. One has $q(F)=30$  {\rm \cite[Prop.~3.6]{CoMuWe}, \cite[Prop.~1]{Let}, \cite{Te}};
\item[$g=3$:]  $X\to Q$ is the general double cover of a smooth quadric $Q\subset\PP^4$ branched along a smooth surface $S\subset Q$ of degree $8$. One has $q(F)=30$  {\rm \cite{Kur}};
\item[$g=4$:]  $X$  is the general complete intersection of a quadric and a cubic hypersurfaces in $\PP^5$. One has $q(F)=20$, $c_1(F)^2 = 23355$, $c_2(F) = 11961$, $p_a(F) = 2942$ \cite[Cor.~18, Thm.~20]{IM2}, \cite{Mar};
\item[$g=5$:] $X$ is the general complete intersection of three quadrics in
$\PP^6$ {\rm \cite{Bea, PB, We}. One has $q(F)=14$, and $X$ is uniquely determined by $F$ \cite{Bea0, De-2, FS2, Las},  \cite[\S 5]{T4}}.
\end{itemize}
\item[$g=6$:]  $F$ is irreducible for any smooth $X=\rm{Gr}(2;5)\cap\PP^7\cap Q\subset\PP^9$ \cite[Cor.~8.3]{DIM1}.  For the general such $X$,  the Abel-Jacobi mapping $A(F)\to J(X)$  is an isogeny,  $q(F)=10$, the image of $F$ in $J(X)$ is algebraically equivalent to $2\Theta^8/8!$ where $\Theta\subset J(X)$ is a Poincar\'e divisor, and $X$ is uniquely determined by $F$ \cite[Thms.\ 0.15, 0.16, 0.18]{Lo2}, \cite{Ili92};
\item Consider a special Gushel threefold $\pi\colon X\to Y$ of genus $6$. Suppose $X$ is general, that is, $X$ is the double cover of a quintic threefold
$Y=\Gr(2;5)\cap\PP^6$ in $\PP^9$  branched along a quadric section  $Y\cap Q$, where $\PP^6\subset\PP^9$ and $Q$ are general. 
Then one has $F=F'\cup F''$ where  $F'=\pi^*(F_1(Y))$ is rational and $F''$ is a non-normal irreducible surface  with
 $q(F'')=10$. For the normalization $\tilde F''$ of $F''$, the Abel-Jacobi mapping $A(\tilde F'')\to J(X)$ is an isomorphism {\rm \cite{Ili94a, Ili94}}.
\item[$g\ge 7$:] $F$ is a smooth irreducible surface with $q(F)=h^{2,1}(X)$, cf.\ Theorem {\rm \ref{thm:q-conics}} and \cite[12.2]{IP}. 
\end{itemize}\end{itemize}
\ethm

Similar results were established in \cite{DIM0, DIM2, FlSe, GLN, MT1, Rad, Ti4,Ti5,  T3, Voi5, We2} for various, possibly singular, Fano threefolds with $\rho=1$ and for various families of curves.
See also, e.g., \cite{AvBVa2, BBR, BlMu, Co, Do0, Do, GrHa, Kan, KuMaMa, Las, MaSt, MaSt1, Me, Merin, OG, PB, PrZa, T4, Voi8}, \cite[Thm.~8.2.1]{IP} and Theorem \ref{thm: cylinder-map-ci} below for some variations and higher dimensional analogs. The classification results for Fano threefolds with $\rho\ge 2$ can be found in \cite{MoMu1}; a part of this classification for $\rho=2$ is compressed in \cite{Ili93} using the Mori theory. \\

The classical Torelli theorem says that any smooth projective curve is uniquely determined by its polarized Jacobian variety; see, e.g., \cite[Ch.~5, Thm.~4.1]{Hu3}. There are different Torelli type theorems for higher dimensional varieties. For instance, for the smooth cubic hypersurfaces of dimension at least three the following holds.

 \bthm\label{thm:Torelli} {\rm (\cite{BD, Cha, CG, Has0, Loo, T3, Voi8}, \cite[Ch.~3, Prop.~2.10, Ch.~3, Prop.~2.10 and Thm.~4.3]{Hu3})}. Assume $X, X' \subset \PP^n$, $n\ge 4$, are smooth cubic hypersurfaces, and let $F=F_1(X)$ and $F'=F_1(X')$ be their Fano varieties of lines endowed with the natural Pl\"{u}cker polarizations $\mathcal{O}_F (1)$ and $\mathcal{O}_{F'} (1)$, respectively.
Then $X \cong X'$ if and only if $(F, \mathcal{O}_F (1)) \cong (F', \mathcal{O}_{F′} (1))$ as polarized varieties. For $n \neq 5$ this is equivalent to $F \cong F'$ as unpolarized varieties.
\ethm

\subsection{The cylinder homomorphism}
The results of the previous subsection are ultimately related to the studies of various cylinder homomorphisms.

\bdefi\label{def: cylinder-map}
Let $X\subset \PP^n$ be a projective variety of dimension $m$, and let $\pi\colon\mathcal{C}\to S$ be a family of irreducible curves in $X$ over an irreducible base $S$. Then the \emph{cylinder homomorphism} associated with $\mathcal{C}$ is defined as follows:
\[\Psi_{\mathcal{C}}\colon H_{m-2}(S,\ZZ)\to H_m(X,\ZZ),\qquad \gamma\mapsto \bigcup_{s\in\gamma} C_s\,,\]
where $\gamma$ is a topological $(m-2)$-cycle in $S$ and $C_s=\pi^{-1}(s)$. We let $\Psi_{\mathcal{C}, \QQ}\colon H_{m-2}(S,\QQ)\to H_m(X,\QQ)$ be the induced homomorphism.
\edefi

For instance, if $S=F_1(X)$ then $\pi\colon\mathcal{C}\to S$ is the universal family of lines in $X$, and if $S=F_c(X)$ then $\pi\colon\mathcal{C}\to S$ is the universal family of conics in $X$, etc. Choosing for $S$ the Fano scheme of lines $F_1(X)$ and letting $\Psi_1, \Psi_{1,\QQ}$ be the associated cylinder homomorphisms, we have the following. 

\bthm\la{thm: cylinder-map-hypersurf} {\rm (\cite{Shim})} Let $X$ be a hypersurface in $\PP^n$ of degree $d$. Assume $n\ge 4$, $d\le n-1$, and either $X$ is general, or $d=3$ and $X$ is smooth. Then the following holds.
\begin{itemize}
\item[$\text{Even}\,\,n$:]
$\Psi_1$ is an isomorphism modulo torsion;
\item[$\text{Odd}\,\,n$:] $\dim(\ker(\Psi_{1,\QQ}))\le (n-3)/4$, and $\Psi_1$ is surjective for $d\le (n+5)/2$.
\end{itemize}
In particular {\rm \cite{BD, CG}}, 
\begin{itemize}
\item[$d\le 4$:] $\Psi_1$ is an isomorphism modulo torsion for $n=5$;
\item[$d=3$:] $\Psi_1$ is surjective, and $\Psi_1$ is an isomorphism modulo torsion for even $n$;
\item[$d=3$:] $\Psi_1$ is  an isomorphism for $n=4,5$.
\end{itemize}
\ethm

Recall \cite[Sect.~1]{Manin1}, \cite{Shim1} that a smooth complete intersection $X\subset \PP^n$ of multidegree $(d_1,\ldots,d_s)$ is a Fano variety of index $\iota$ if and only if $\iota:=n+1-\sum_{i=1}^{s} d_{i}>0$. 

\bthm\la{thm:conics-in-Fci} 
Let $X\subset\PP^n$ be a smooth Fano complete intersection. 
Then the following holds.
\begin{itemize}\item  $X$ is covered by conics \cite{CoMu, CoMu1, Lew}; 
\item$X$ is covered by lines provided $\iota(X)\ge 2$ \cite[Lect.\ 4, Proof of Lemma 1]{T3};
\item the lines in $X$ sweep out a hypersurface provided $\iota(X)=1$  \cite{Shim}. 
\end{itemize}
\ethm

The following theorem comprise certain results of \cite{Lew1, Lew2} for the case of Fano complete intersections of index 1.  

\bthm\la{thm: cylinder-map-ci}  {\rm (\cite{Shim1})}
Let $X$ be a general Fano complete intersection in $\PP^n$ of dimension $k\ge 3$. Then 
\begin{itemize}\item both $\Psi_{1,\QQ}$ and $\Psi_{c,\QQ}$ are surjective; 
\item
 if $k=2s-1$ is odd, then the Abel-Jacobi mappings \[J^{s-1}(F_1(X))\to J^{s}(X)=J(X)\quad\text{and}\quad J^{s-1}(F_c(X))\to J^{s}(X)=J(X)\] are surjective.
\end{itemize}
\ethm

See, e.g., \cite{BlMu, KLMS, KL, Lew1, Lew2, Lew3, Pan, Shio, Voi4} and the literature therein for the cylinder mappings on cycles of other intermediate dimensions. See also \cite{Bel} on the intermediate jacobians of conic bundles. See, e.g., \cite{DeKu1}--\cite{DeKu4} for the Gushel-Mukai varieties and their  intermediate jacobians.

\section{Counting rational curves }\label{ss:rat-curv}
\subsection{Varieties of rational curves in hypersurfaces}\label{ss:rat-curves-hypers}
Given  a hypersurface $X$ of degree $d$ in $\PP^n$, we let $R_e(X)$ be the space of smooth rational curves of degree $e$ in $\PP^n$ lying in $X$. This is an open subscheme of the Hilbert scheme ${\rm Hilb}^{et+1}(X)$. The number
\begin{equation}\label{eq:expected-dimension} \mu_e=\mu_e(d,n)=(n+1-d)e+n-4\end{equation}
is called the \emph{expected dimension} of $R_e(X)$. Notice that \[\mu_2=\mu_2(d,n)=-\epsilon(d,n)=3n-2d-2\,.\] 

\begin{thm}\label{thm:only-lines} Let $X$ be the general hypersurface of degree $d\ge 2$ in $ \PP^n$, where $n\ge 3$. Then the following holds.
 \begin{itemize} 
\item[{\rm (a)}]  
 $R_e=\emptyset$  if $\mu_e<0$ and $e\le d+1$ {\rm \cite[Thm.\ 1.1]{Fur1}}.
\item[{\rm (b)}]  \begin{itemize} \item $R_e(X)$ is smooth of dimension $\mu_e$ if $\mu_e \ge 0$ and either $e\le 3$, or $2e\le d+3$;  
\item $R_e(X)$ is an integral, locally complete intersection scheme of dimension $\mu_e$
if $d \le n-2$  {\rm \cite{BehKu, RiYa1}, \cite[Thm.\ 1.1]{Fur1}, \cite[Thm.\ 1.1]{HRS1}}, see also  \cite{BrVi, HaSt, St0, St}; 
\item $R_e(X)$ is  irreducible, generically smooth and of dimension $\mu_e$ if $n \ge 4$,  $d \le n-1$, and $e \le d-1$ \cite{HRS1}, \cite[Thm.~26]{ZRa3}, \cite{Ts}.
\end{itemize}
\item[{\rm (c)}]  
If $2d \le n+1$ then through any point of $X$ passes a family of degree $e$ rational curves of dimension $e(n+1 - d) - 2\ge ed$. 
In particular, through any point  of $X$ passes a $2(n-d)$-dimensional family of smooth conics {\rm \cite{HRS1}}. 
\item[{\rm (d)}] Let $R$ be a sweeping component of $R_e(X)$, 
that is, the corresponding rational curves sweep out an open subset of $X$. If $(n + 1)/2 \le d \le n -3$, then $R$ is not uniruled  \cite[Thm.~1.1]{Beh3}.
\item[{\rm (e)}] 
For any $n \ge 4$, $d \le n$, and $e \le 2n-2$ there exists on $X$ a rational curve of degree $e$ with balanced normal bundle  {\rm \cite{CoRi2}, \cite[Thm.~24]{ZRa3}}.
\end{itemize}
\end{thm}

Notice that the dimension of $R_e(X)$ can be strictly larger than $\mu_e$ for
particular smooth hypersurfaces $X$. For instance \cite{Ver}, the family of lines
on the Fermat quartic in $\PP^4$ is two-dimensional, while the general quartic in $\PP^4$
carries a one-parameter family of lines. See further examples in \cite[Ex.\ 3.17--3.18]{Fur2}. See also \cite{BehKu, CoSt} 
 for results on the Gromov-Witten invariants.

Concerning arbitrary smooth hypersurfaces, we have the following results.

\begin{thm}\label{thm:smooth}  Let $X$ be a hypersurface  of degree $d$ in $\PP^n$.
\begin{itemize}
\item[{\rm (a)}] 
 If $X$ is smooth along $C$ for some $C \in R_e(X)$ then $\dim_C R_e(X)\ge\mu_e=\chi(\mathcal N_{C/X})$ {\rm \cite[Rem.\ 3.2]{Fur2}, \cite[II, Thm.\ 1.2]{Kol}}.
\item[{\rm (b)}]
 For any smooth hypersurface $X\subset \PP^n$ of degree $d$ where either $d=3$ and $n\ge 5$, or $d\ge 4$ and $n\ge 2^{d-1}(5d-4)$, 
  the scheme $R_e(X)$ is irreducible of the expected dimension $\mu_e$ {\rm \cite{BrVi, CoSt}}. 
 \item[{\rm (c)}]  Assume  $X$ is smooth and $d+2k - n \ge 3$. Then the quadrics of dimension $k$ sweep out a subvariety of dimension at most $n-k-1$ in $X$  {\rm \cite{Beh0}}. \footnote{In particular, for $d\ge n+1$ the conics contained in $X$ do not cover $X$. The latter  can be shown directly.}
 \item[{\rm (d)}] Assume $d=n\ge 5$. Then there is a countable set of closed, codimension two subvarieties of $X$  such that the
image of any generically finite, regular morphism from a del Pezzo surface to $X$ is contained in one of these \cite[Thm.~1.4]{BehSt}.
\end{itemize}
\end{thm}

In \cite{St} one can find restrictions under which the Kontsevich moduli space $\overline{\mathcal{M}}_{0,0}(X, e)$ is of general type. See also \cite{Pan2} for studies on the Kontsevich moduli spaces of rational curves through marked points in a Fano complete intersection variety. 

  Recall that a projective variety $X$ is \emph{rationally connected} if each pair of closed
points of $X$ is contained in a rational curve, see  \cite{De0, Kol}. 
Notice that any Fano variety $X$ is rationally connected and covered by the rational curves of degree $\le \dim(X)+1$. These curves generate the Mori cone of effective $1$-cycles on $X$ \cite{Cam, KMM0, Mo1}, \cite[Sect.~IV.3, Cor.~IV.1.15]{Kol}. See also \cite{dJSt} for criteria of \emph{simple rational connectedness}.
 
The following results concern rational curves and rational surfaces in smooth complete intersection Fano varieties. 
It is known that a general such variety of sufficiently small multidegree is unirational \cite{AM, Kol, Kol1, KMM, PaSr, Wal}.

\bthm\la{rat-simple-conn} Let $X\subset\PP^n$ be a smooth complete intersection  of
type $\underline{d}=(d_1,\ldots,d_s)$. Then the following holds.
\begin{itemize} 
\item[{\rm (a)}]  $X$ is rationally connected if and only if $\sum_{i=1}^s d_i\le n$, that is, $X$ is a Fano variety  {\rm \cite{dJSt}}.
\item[{\rm (b)}]  Assume $\omega_X\simeq \mathcal{O}_X (-1)$ that is, $\sum_{i=1}^s d_i= n$.
 Let $S$ be a smooth variety of dimension $2 \le \dim(S) \le \dim(X) -2$ with $\omega_S^\vee$ nef. Consider a generically finite morphism $f \colon S \to X$. In the case $\omega_S\simeq \mathcal{O}_S$ assume further that $f (S)\subset\PP^n$ is linearly non-degenerate.  Suppose $f$ extends to a morphism $F\colon \mathcal{S}\to X$, where $\mathcal{S}\to \mathcal{B}$ is a deformation family containing $S$ as a fiber. 
Then the image $F(\mathcal{S})$ is contained in a subvariety of codimension at least two in $X$ \cite[Thm.~1.1]{Tes08}.
\end{itemize}
\ethm

Notice that statement (b) generalizes Theorem \ref{thm:smooth}(d).  See also \cite{Va, Va1} for count of rational and elliptic curves on rational surfaces and in projective spaces, \cite{Sh0, Sh1} for count of rational curves in Fano threefolds and Fermat hypersurfaces, and \cite{Zah} for count of rational curves in del Pezzo manifolds. 

\subsection{Twisted cubics in complete intersections}\la{ss:twisted-cubics}
The next results concern enumeration of twisted cubics on Fano complete intersections. 

\bthm\la{thm:twisted-cubics-in-c.i.} \begin{itemize}\item[{\rm (a)}] 
Let $\Sigma_{tc}(d,3)\subset\Sigma(d,3)$ be the locus of degree $d$ surfaces in $\PP^3$ which contain twisted cubic curves. 
There is an explicit expression of the degree of $\Sigma_{tc}(d,3)$ as a polynomial in $d$ of degree $24$ {\rm \cite[\S 8]{MaRoXaVa}}. 
\item[{\rm (b)}]
Let $X\subset\PP^{n}$ be a smooth complete intersection of type $(d_1,\ldots,d_s)$ where
$$ \sum_{i=1}^s (d_i-1)=\frac{n-s}{3}+1\,.$$ Then the number of twisted cubics in $X$ passing through three general points in $X$ equals
$$ \prod_{i=1}^s ((d_i-1)!)^2d_i!\,.$$
In particular, through three general points of a smooth cubic threefold $X\subset\PP^4$ pass exactly $24$ twisted cubics in $X$  {\rm \cite{Bea1}}. 
\end{itemize}
\ethm

See also \cite{PS}, the survey article \cite{Pi}, and the references therein. See, e.g., \cite{BCS, dJSt, Has1, Hu4, Hu3, KLSV, Kuz1, LaLeMa, Lehn, LLSVS, LPZ, Sac, ShYi} for twisted cubics on cubic fourfolds and related hyper-K\"{a}hler varieties. 

\section{Hypersurfaces with few rational curves}\label{sec:few-curves}
\subsection{Rational curves on ${\rm K}3$ surfaces}\label{ss:rat-curv-K3} 

Recall the following theorem.

\bthm\la{thm:genus-g-curves-in-quartic-surf} {\rm (\cite{Mo})} Let $d>0$ and $g \ge 0$ be integers. There is a smooth curve $C$ of degree $d$ and genus $g$ lying in a smooth quartic surface $X$ in $\PP^3$ if and only if either $g  =  d^2/8 +1$, or $g < d^2/8$ and  $(d,g)\neq  (5,3)$.
\ethm

For surfaces in $\PP^3$ containing smooth elliptic quartic curves, there are the following facts. 

\bthm\la{thm:ell-quartic-curves-in-surfaces}  Let $\Sigma_{\rm ell,4}(d,3)$ be the  locus of surfaces of degree $d$ in $\PP^3$
containing an elliptic quartic curve. Then the following holds.
\begin{itemize}\item[{\rm (a)}]
$\Sigma_{\rm ell,4}(4,3)$  is a hypersurface in $\Sigma(4,3)$
of degree $38475$  {\rm \cite[\S 3.2]{CuLoVa}}. 
\item[{\rm (b)}]  There is an explicit  expression of $\deg(\Sigma_{\rm ell,4}(d,3))$ for $d\ge 5$  as a polynomial in $d$ of degree $32$ obtained via Bott's residue formula {\rm \cite[\S 4.3]{CuLoVa}}. 
\end{itemize}
\ethm

Recall (see, e.g., \cite{Lew}) that a smooth hypersurface
$X$ in $\PP^n$ of degree $d\ge n+1$ cannot be covered by rational curves. 
This concerns, in particular, smooth quartic surfaces in $\PP^3$, and holds, 
more generally, for any ${\rm K}3$ surface. Nonetheless, a projective ${\rm K}3$ surface carries infinitely many rational curves, 
see (a) in the next theorem. This generalizes a previous partial result due to Bogomolov, 
Mumford, Mori-Mukai \cite{MoMu}, see also \cite{BaMK, BoHaTs, BoTs1, BoTs2, Ch1, Has, LiLi, DMOR}. 

\begin{thm}\label{thm:rat-curves-on-K3} Let $X$ be a projective ${\rm K}3$ surface over an algebraically closed field.  Then the following holds.
\begin{itemize}\item[{\rm (a)}]  
 $X$ contains infinitely many rational curves {\rm \cite{CGL, Tay}; cf. also \cite{CGL1}}.
\item[{\rm (b)}] 
 Consider 
the subset $\mathcal{S}_g$ of the moduli scheme of the $K3$ surfaces of genus $g$ which parameterizes the surfaces 
$X$ such that the union of rational curves 
on $X$ is dense in the Hausdorff topology. Then $\mathcal{S}_g$  is of the second Baire category {\rm \cite{ChLe}}.  
\item[{\rm (c)}]  
Any rational curve on a general  ${\rm K}3$ surface $X$ of genus $g\ge 2$ is nodal.
 \item[{\rm (d)}]  Given $d\in\NN$, for any  $h>84d^2$ and for any  $K3$ surface $X$ of degree $h$
over a field $k$ of characteristic $p\neq 2,3$ the number of rational curves in $X$ of degree at most $d$ does not exceed  $24$. This upper bound is exact for any $d\ge 3$ \cite{Mi, RaSc5}. 
\end{itemize}
\end{thm}

Similar results hold also for the Enriques surfaces, see \cite{RaSc7} and the references therein.

See \cite{Hu1, Hu2} for a discussion related to (b).  The generality assumption in (c) is essential; indeed, there are smooth quartic surfaces in $\PP^3$ which contain the rational 3-cuspidal plane quartic {\rm \cite{Ch1, Ch2}}.
As for the count of curves on a ${\rm K}3$ surface in terms of the Gromov-Witten invariant, see, e.g., \cite{Bea2, LeLe, MaPa,  Wu, YaZa}.

Notice that there are complex analytic $K3$ surfaces with no algebraic curve. 
However, any $K3$ surface $S$ carries a transcendental entire curve, that is, the image of a nonconstant holomorphic map $\C\to S$ \cite{Can}.

\subsection{Rational curves in hypersurfaces}\label{ss:few-rat-curves-hypers}
The next theorem deals with hypersurfaces which contain few rational curves. 

\begin{thm}\label{thm:few-rat-curves} Let $X$ be a hypersurface  of degree $d$ in $\PP^n$. Then the following holds. 
\begin{itemize}
\item[{\rm (a)}] 
For $n\ge 3$ and $2(d+1)\ge 3(n+1)$ any rational curve in a very general $X$ is contained in the maximal subvariety $L(X)$ of $X$ swept out by lines in $X$ {\rm \cite[Cor.\ 3.2]{RiYa2}}.
\item[{\rm (b)}] In the range  $n\ge 7$ and $d\in [\frac{3n+1}{2},2n-3]$ the general $X$ contains lines but no other rational curves   {\rm \cite[Thm.\ 1.3]{RiYa2}}. This is true as well for $n=6$ and $d=2n-3$ \cite{Pa1}, but fails for  
a general quintic threefold in $\PP^4$. 
\item[{\rm (c)}] For $n=5$ and $d\ge 2n-3$ the general $X$ contains just a finite number of rational curves of any given degree (that is, $X$ does not contain any one-parameter family of rational curves) {\rm \cite{Vo1}}. 
\item[{\rm (d)}]  If $X$ is very general, $n \ge 4$, and $d \ge 2n -2$ then $X$ contains no rational curve {\rm \cite{Vo1}}. 
\end{itemize}
\end{thm}

Notice that (d) fails for $n=3$.  Indeed,  by  Theorem \ref{thm:rat-curves-on-K3}(a),
any smooth quartic surface in $\PP^3$ contains a sequence of rational 
curves of growing degrees. 
Statement (c) is a strengthening of the previous results in \cite{Cl, Ei, Ei2}. The validity of an analog of (c) for $n=4$ and $e\ge 13$ is still open; 
this  is the famous Clemens Conjecture. 

\subsection{Clemens' Conjecture} \label{ss:high-degree-curves} 
This conjecture \cite{Cl} suggests that a general quintic threefold in $\PP^4$ 
contains a positive finite number of smooth rational curves of any given degree, 
and the scheme of such curves is reduced. 
The next theorem gives a brief summary of some results on Clemens' Conjecture. 

\begin{thm}\label{thm:Clemens-Conj}  \begin{itemize}
\item[{\rm (a)}] 
The Clemens Conjecture holds for curves of degree $\le 12$ {\rm \cite{BaFo, Cot1, Cot2, JK1, Ka}}. 
\item[{\rm (b)}]  For any $d\ge 1$ a general quintic threefold in $\PP^4$ 
contains a smooth rational curve $C$ of degree $d$ {\rm \cite{Ka}}. 
\item[{\rm (c)}]  Any smooth rational curve $C$ in a general quintic threefold $X$ in $\PP^4$ is embedded with normal bundle $\mathcal O_{\PP^1} (-1) \oplus\mathcal O_{\PP^1} (-1)$. Any singular rational curve in $X$ is a plane $6$-nodal quintic {\rm \cite{Cot1, Cot2, JK1, Ka, Ni}}.
\item[{\rm (d)}] 
The number of smooth rational curves of degree $1,2,3,4,...10$ in a general quintic threefold in $\PP^4$ is, respectively,
$$2875, \,\,\, 609250, \,\,\, 317206375, \,\,\, 242467530000,\ldots,704288164978454686113488249750\,,$$
where the number $2875$ of lines is due to Schubert \cite{Sch}; see \cite{Pi}, \cite[\S 10.6]{BePe} and the references therein.
\end{itemize}
\end{thm}

For instance \cite{Ka}, \cite{Shin2}, a general hypersurface of degree $d > \frac{3}{2} n - 1$ in $\PP^n$ does not contain any smooth conic; however, a general quintic threefold in $\PP^4$ does.  

 For any natural number $d\ge 1$ there is a Mirror Symmetry prediction for the number of smooth rational curves of degree $d$ in a general quintic threefold in $\PP^4$. Actually, these virtual numbers count pseudoholomorphic curves in a general almost complex symplectic deformation of the quintic threefold via  quantum cohomology, see, e.g., \cite{BePe, Kon, KoMa, McDSa}. 
 
The same methods work for certain smooth Calabi-Yau complete intersections (CICY, for short). Besides the quintic threefolds in $\PP^4$, there are exactly 4 types of  smooth CICY threefolds of type, respectively, $(3,3)$ and $(2,4)$ in $\PP^5$, $(2,2,3)$  in $\PP^6$, and $(2,2,2,2)$ in $\PP^7$. The Mirror Symmetry prediction for the number of smooth rational curves of degree $d\le 10$ in a general CICY threefold can be found in \cite{LiTe}. For $d\le 6$ this prediction gives the correct number  of curves, see \cite{H3,Li}. 
The next theorem addresses rational and elliptic curves in general Calabi-Yau complete intersection threefolds. 

\bthm\la{thm:ell-curves-in-CICY-3folds} {\rm (\cite{EkJoSo, Kl, Kn2})}  
Let $X$ be a general CICY threefold.
Then for $g=0,1,2,3$ there is an integer $d_g\ge 0$, where $d_0=0$, such that, for any $d > d_g$,  $X$ contains an isolated smooth curve of degree $d$ and genus $g$.  
\ethm

 In \cite{Kl}, some of these results are claimed to hold for any smooth CICY threefolds. However, there is a gap in the proof  in \cite{Kl}; see \cite{Kn2}.
Similar facts hold for certain higher genera curves under more severe restrictions, see  \cite[Thm.\ 1.2]{Kn2}. 

The following theorem gives a short account of sporadic results for curves in hypersurfaces. 

\bthm \begin{itemize}
\item A  very general  hypersurface of
degree $d \geqslant 2n - 1$ in $\PP^n$ does not contain any smooth elliptic curve  {\rm \cite{Wang1}; cf.\ \cite{Clm1}}.
\item  The degree of an elliptic curve on a very general heptic hypersurface $X$ in $\PP^4$ is a multiple of $7$ {\rm \cite{FR}}. 
\item Let $X$ be the general heptic hypersurface in $\PP^5$. Then $X$ does not contain any rational curve of degree $d\in\{2,\ldots,16\}$  {\rm \cite{Cot3, HJ, Shin1}}, any smooth elliptic curve of degree $e \le 14$, and any smooth curve $C$ of degree $e \le16$ and genus $1 \le g \le 3$ 
provided the dimension of the linear span of $C$ is not equal to $3$ \cite{Ba}. 
\item  A general hypersurface of degree $54$
in $\PP^{30}$ does not contain any rational quartic curve {\rm \cite{Wang2}}.
\end{itemize}
\ethm

The second statement goes in the direction of the conjecture of Griffiths and Harris \cite{GrHa1} which says that for a very general hypersurface of degree $d \ge 6$ in $\PP^n$ and for any curve $C$ in $H$ one has $d|\deg(C)$. This is true for $n=3$ due to the Noether-Lefschetz theorem; see \cite{BaCaCi} for further results. See also the survey article \cite{He} on the role of the Calabi-Yau,  in particular, CICY varieties in physics.

 \section{Counting curves of higher genera and hyperbolicity}\label{ss:higher-degree-rat-curves}
It is worthwhile to compare previous results with the following finiteness theorems related to Kobayashi hyperbolicity 
and the Green-Griffiths-Lang Conjecture (see below). 
  
 \bthm\la{thm:Kob-hyp} \begin{itemize}\item[{\rm (a)}]  
Consider a projective variety $X\subset\PP^n$. If $X$ is Kobayashi hyperbolic then
 there exists $\varepsilon>0$ such that for any curve $C$ of geometric genus $g$ in $X$ one has  {\rm \cite{Dem1}}
\be\label{eq:alg-hyp}
2g-2 \ge \varepsilon \deg(C)\,.\ee
Consequently, the curves of a given geometric genus in $X$ form a bounded family. 
\item[{\rm (b)}] A general hypersurface $X$ in 
$\PP^n$ of degree $d\ge 16(2n-3)^5(10n-11)$ is Kobayashi hyperbolic  {\rm \cite{BeKi, Bro, Dem2, YDe, Mer, RiYa3, Si}}. 
The latter holds for $n=3$ starting with $d=18$ \cite{Pau} and for $n=4$ starting with $d=593$ \cite{DT}. 
 \end{itemize}
\ethm

A weaker form of \eqref{eq:alg-hyp} called \emph{algebraic hyperbolicity} implies the absence of rational and elliptic curves; see, e.g., \cite{Ch0, CL, CoRi, HaIl}. 
See also \cite{BD, Dar} for logarithmic versions of (b). There are examples of smooth hyperbolic surfaces in $\PP^3$ of any given degree $d\ge 6$  \cite{CZ0, Du, SZ1, SZ3, SZ4, Za}  and of hyperbolic hypersurfaces in $\PP^n$ of degree $d\sim n^2/4$ \cite{Huy, SZ2}.  

 In the direction of algebraic hyperbolicity of general hypersurfaces, the following holds.

\bthm\label{thm:gen-type} {\rm (\cite{Pa2})} For $n\ge 6$ and for a very general hypersurface $X$ in $\PP^n$ of degree $d\ge 2n-2$, any subvariety $Y\subset X$ is of general type.
\ethm

Cf.\ also \cite{Xu, Xu4, Xu5}. 
  
\bthm\la{thm:BM} \begin{itemize}\item[{\rm (a)}]  
The number of rational and elliptic curves on 
a minimal smooth surface of general type with $c^2_1 > c_2$ is 
bounded above by an effective function of $c_1$ and $c_2$ {\rm \cite{Bog, Des, Mi}, \cite[Thm.~10.1]{La2}}.
\item[{\rm (b)}] Let $X$ be a very general surface  of degree $d \ge 5$ in $\PP^3$. Then
for any curve $C$ of geometric genus $g$ on $X$ one has \be\la{eq:CR-bound} 2g-2\ge \max\Big\{d^2-3d-6,\,\left(d+\frac{1}{d}-5\right)\deg(C)\Big\}
\,.\ee In particular, for $d=5$ one has 
$$2g-2 \ge  \max\Big\{4,\frac{1}{5} \deg(C)\Big\}\,,$$ which yields $g\ge 3$.
The absolute lower bound $2g-2\ge d^2-3d-6$ is sharp; for $d \ge 6$ it is achieved by the tritangent hyperplane sections only  
{\rm \cite{CL, CoRi, Xu}}. 
\end{itemize}
\ethm

Cf.\ also \cite{LuMi, LM1, Wu} for (a), \cite{CoRi1, HaIl} for (b). 

For elliptic curves on general $K3$ surfaces and Fano schemes of lines on general cubic fourfolds, the following hold.

\bthm\label{thm:ell-curves} {\rm (\cite{NeOb})}
\begin{itemize}\item[{\rm (a)}] Let $X$ be a very general $K3$ surface with primitive curve class $\beta\in H^2(X, \ZZ)$ of self-intersection $\langle\beta, \beta\rangle=2h-2$, $h\in\ZZ_{\ge 0}$.
Then the moduli space of elliptic curves on $X$ from the class $\beta$ is a smooth curve. The number $n_{\beta,j}$ of such curves with fixed general $j$-invariant depends only on $h$. It admits an expression in terms of two particular Gromov--Witten invariants, which can be computed explicitly for any given value of $h$. 
\item[{\rm (b)}] A general Fano variety of lines on a cubic fourfold in $\PP^5$ contains precisely $3780$ elliptic curves of minimal degree and of fixed general $j$-invariant.
\end{itemize}
\ethm

Notice that there are exactly $6383765416$ elliptic quartics meeting $16$ general lines in $\PP^3$ \cite{AvBVa1}. 

For higher genera curves on $K3$ surfaces, we have the following theorem. 

\bthm\label{thm:K3-g}  {\rm\cite[Cor.~2]{Nis}}
A generic projective $K3$ surface contains infinitely many $g$-dimensional families of irreducible immersed curves of geometric genus $g$, for any positive integer $g$.
\ethm

For the genera of curves on smooth surfaces in $\PP^3$ the following is known.

\bthm\label{thm:gaps} For $d\geqslant 4$ let ${\rm Gaps}(d)$ be  the set of all the non--negative integers
which cannot be realized as geometric genera of irreducible curves
on a very general surface of degree $d$ in $\PP^ 3$. Then ${\rm Gaps}(d)$ is the union of
finitely many disjoint and separated integer intervals ${\rm
Gaps}_j(d)$, $j=0,1,\ldots$. One has:
\begin{itemize} 
\item ${\rm
Gaps}(5)= \{0,1,2\}$  {\rm \cite{Xu}}; \item
${\rm
Gaps}_0(d)=\left[0, \; \frac{d(d-3)}{2} - 3\right]$ for all $d\geqslant 5$  {\rm \cite{Xu}}; 
\item
${\rm Gaps}_1(d)=\left[\frac{d^2-3d+4}{2}, \; d^2
- 2d - 9\right]$ for all $d\geqslant 6$   {\rm \cite{CiFlZa1}}.  
\end{itemize}
\ethm

In the other direction, for the existence of subvarieties with a given geometric genus, we have the following result. 

\begin {thm}\label{thm:main} {\rm (\cite{ACH} for $s=1$, \cite[Thm.\ 0.1]{CiFlZa0})} Let $X$ be a smooth  irreducible projective variety of dimension $n>1$, let $L$ be a very ample divisor on $X$, and let $s\in\{1,\ldots,n-1\}$. 
Then there is an  integer $p_{X,L,s}$ (depending on $X$, $L$ and $s$) such that 
for any $p\geqslant p_{X,L,s}$ one can find an irreducible subvariety $Y$ of $X$ of dimension $s$ with at most ordinary points of multiplicity $s+1$ as singularities such that $p_g(Y)=p$. Moreover, one can choose $Y$ to be a complete intersection $Y=D_1\cap\ldots\cap D_{n-s}$, where $D_i\in |L|$ for $i=1,\ldots,n-s-1$ are smooth and transversal and $D_{n-s}\in |mL|$ for some $m\geqslant1$ is such that $Y$ has  ordinary singularities of multiplicity $s+1$. 
\end{thm} 

\bcor\label{cor:gaps}  {\rm (a)}  For any integer $d\ge 4$, there exists an integer $c(d)$ such that, for any smooth surface $S$ in $\PP^3$ of degree $d$ and any integer $g\ge c(d)$, $S$ carries a reduced, irreducible nodal curve of geometric genus $g$, whose nodes can be prescribed generically on $S$ {\rm \cite{ACH}, \cite[Cor.\ 3.1]{CiFlZa0}}.

{\rm (b)} For any positive integer $d$ and for any non-negative integer $g$, there is a smooth surface $S$ in $\PP^3$ of degree $d$ and an irreducible, nodal curve $C$ on $S$ with geometric genus $g$  {\rm \cite[Thm.\ 3.3]{CiFlZa0}}.\ecor

Recall the famous

\smallskip

{\bf Green-Griffiths-Lang Conjecture.}   {\rm(\cite{GG, La})}
\emph{Let $X$ be a projective variety of general type. Then there exists a proper closed subset $Y\subset X$ which contains any subvariety $Z\subset X$ not of general type and the image of any nonconstant entire curve $\CC\to X$.}

\smallskip

The conjecture is fixed in the particular case of general projective hypersurfaces.

\bthm\label{thm:GGL-conj} {\rm(\cite{Berc, BeKi, BJK, Dem2, DMR, Si})}
For a generic projective hypersurface $X\subset\PP^n$, $n\ge 2$, of degree $d\ge 16(n-1)^5(5n-1)$ there exists a proper subvariety $Y\subset X$ which contains the image of any nonconstant entire curve $\CC\to X$, hence also any rational or elliptic curve on $X$ and the images of abelian varieties.
\ethm

There are the following related conjectures.

\smallskip
 
{\bf Conjecture.} (C.~Ciliberto, F.~Flamini, and M.~Zaidenberg \cite{CiFlZa}) \emph{There exists a strictly growing function $\varphi\colon\NN\to\NN$ such that the number of curves of geometric genus $g\leqslant\varphi(d)$ in any smooth surface $S$ of degree $d\geqslant 5$ in $\PP^3$ is finite and bounded by a function of $d$. }
 
\smallskip
 
 {\bf Conjecture.} {\rm (C.~Voisin \cite{Voi3})} \emph{ Let $X\subset\PP^n$ be a very general hypersurface of degree $d\geqslant n+2$. Then the degrees of rational curves in $X$ are bounded.}
 
 \smallskip
 
 {\bf  Conjecture.}  {\rm (P.~Autissier, A.~Chambert-Loir, and C.~Gasbarri \cite{ACLG})} \emph{ Let $X$ be a smooth projective variety of general type with canonical line bundle $K_X$. Then there exist real numbers $A$ and $B$, and a proper Zariski closed subset $Z \subset X$ such that for any curve $C$ of geometric genus $g$ in $X$ not contained in $ Z$, one has
$\deg_C(K_X) \leqslant A(2g - 2) + B$.}

\smallskip



\end{document}